\documentclass{amsart}

\usepackage{amsmath}
\usepackage{amssymb}
\usepackage{amsthm}
\usepackage[all]{xypic}
\usepackage{verbatim}

\numberwithin{equation}{section}
\theoremstyle{plain}
\newtheorem{theorem}[equation]{Theorem}
\newtheorem{proposition}[equation]{Proposition}
\newtheorem{corollary}[equation]{Corollary}
\newtheorem{lemma}[equation]{Lemma}

\theoremstyle{definition}
\newtheorem{defn}[equation]{Definition}

\theoremstyle{remark}

\newcommand{\beq}{\begin{equation}}
\newcommand{\eeq}{\end{equation}}

\newcommand{\blank}{\mbox{$\underline{\makebox[10pt]{}}$}}

\newcommand{\st}{\left\vert\right.}

\newcommand{\bbar}[1]{\overline{#1}}

\DeclareMathOperator{\Hom}{{Hom}}

\DeclareMathOperator{\Ext}{{Ext}}

\DeclareMathOperator{\Tor}{Tor}
\DeclareMathOperator{\Aut}{{Aut}}

\DeclareMathOperator{\reg}{reg}

\DeclareMathOperator{\im}{Im}

\DeclareMathOperator{\gr}{gr}

\DeclareMathOperator{\sing}{sing}
\DeclareMathOperator{\Cosupp}{Cosupp}

\newcommand{\sh}{\mathcal}
\newcommand{\mc}{\mathcal}
\newcommand{\mf}{\mathfrak}

\DeclareMathOperator{\shHom}{\mathcal{H}\!\mathit{om}}
\DeclareMathOperator{\shExt}{\mathcal{E}\!\mathit{xt}}
\DeclareMathOperator{\shTor}{\mathcal{T}\!\mathit{or}}
\renewcommand{\Lsh}{\mathcal{L}}

\newcommand{\struct}{\mathcal{O}}
 \newcommand{\kk}{{\Bbbk}}
\newcommand{\sA}{\mc{A}}
\newcommand{\sB}{\mc{B}}
\newcommand{\sC}{\mc{C}}
\newcommand{\sD}{\mc{D}}
\newcommand{\sE}{\mc{E}}
\newcommand{\sF}{\mc{F}}

\newcommand{\sH}{\mc{H}}
\newcommand{\sI}{\mc{I}}
\newcommand{\sJ}{\mc{J}}
\newcommand{\sK}{\mc{K}}
\newcommand{\sL}{\mc{L}}

\newcommand{\sO}{\mc{O}}
\newcommand{\sP}{\mc{P}}
\newcommand{\sQ}{\mc{Q}}
\newcommand{\sR}{\mc{R}}
\newcommand{\sS}{\mc{S}}
\newcommand{\sT}{\mc{T}}
\newcommand{\sU}{\mc{U}}

\newcommand{\ZZ}{{\mathbb Z}}

\newcommand{\PP}{{\mathbb P}}
\newcommand{\NN}{{\mathbb N}}

\newcommand{\mb}{\mathbb}

\DeclareMathOperator{\rgr}{gr-\!}
\DeclareMathOperator{\lgr}{\!-gr}

\DeclareMathOperator{\rGr}{Gr-\!}

\DeclareMathOperator{\lmod}{\!-mod}
\DeclareMathOperator{\lMod}{\!-Mod}

\DeclareMathOperator{\lQgr}{\!-Qgr}
\DeclareMathOperator{\lqgr}{\!-qgr}

\DeclareMathOperator{\rQgr}{Qgr-\!}
\DeclareMathOperator{\rqgr}{qgr-\!}

\DeclareMathOperator{\rTors}{Tors-\!}

\DeclareMathOperator{\GT}{GT}

\newcommand{\I}{{\mathbb I}}

\newcommand{\ver}[1]{^{(#1)}}

\DeclareMathOperator{\Supp}{Supp}

\newcommand{\surfdata}{(X, \Lsh, \sigma, \sh{A}, \sh{D}, \sh{C},s)}
\newcommand{\ADCdata}{(X, \Lsh, \sigma, \sh{A}, \sh{D}, \sh{C},s)}
\newcommand{\geomdata}{(X,  \sigma, \Lambda, Z, \Gamma, \Omega)}
\newcommand{\ADCD}{\mathbb{E}}
\newcommand{\surfD}{{\mathbb D}}

\title{Geometric algebras on projective surfaces}
\author{Susan J. Sierra}
\address{Mathematics Department,  Fine Hall, Washington Road, Princeton, NJ 08544}
\email{ssierra@princeton.edu}
\date{\today}

\begin{document}

\begin{abstract}
Let $X$ be a projective surface, let $\sigma\in \Aut(X)$, and let $\Lsh$ be a $\sigma$-ample invertible sheaf on $X$.  We study the properties of a family of subrings,  parameterized by geometric data, of the twisted homogeneous coordinate ring $B(X, \Lsh, \sigma)$; in particular, we find necessary and sufficient conditions for these subrings to be noetherian.  We also study their homological properties, their associated noncommutative projective schemes, and when they are maximal orders.  In the process, we produce new examples of   maximal orders; these are graded and have the property that no Veronese subring is generated in degree 1.  

Our results are used in the companion paper \cite{S-surfclass} to give defining data for a large class of noncommutative projective surfaces. 
\end{abstract}
\maketitle

\tableofcontents

\section{Introduction}\label{INTRO}
Twisted homogeneous coordinate rings and their subalgebras, which we refer to in this paper as {\em geometric algebras}, have been important sources of examples and counterexamples in noncommutative ring theory in recent years; see particularly \cite{AV}, \cite{KRS}, and \cite{S-idealizer}.  The aim of this paper is to analyze the properties of a broader class of geometric algebras than have  been studied so far.  In  \cite{S-surfclass}, the companion paper to this one, we use these algebras to complete an important special case of the classification of noncommutative projective surfaces.

Let us give the geometric data that define the rings under study.  We work over a fixed algebraically closed field $\kk$.  
If $X$ is a projective variety, $\sigma \in \Aut(X)$, and $\Lsh$ is a quasicoherent sheaf on $X$, we will write 
\[\Lsh^{\sigma} := \sigma^* \Lsh.\]

\begin{defn}\label{surfdata}
The tuple $\surfD = \surfdata$ is {\em ADC data} if:
\begin{itemize}
\item $X$ is  a projective surface;
\item $\sigma$ is an automorphism of $X$;
\item $\Lsh$ is an invertible sheaf on $X$;
\item $s$ is a positive integer;
\item $\sh{D}$ is the ideal sheaf of a 0-dimensional subscheme of $X$ such that all points in the cosupport of $\sh{D}$ have distinct infinite $\sigma$-orbits; and
\item $\sh{A}$ and $\sh{C}$ are ideal sheaves on $X$ such that the cosupport of $\sC$ is 0-dimensional and
\beq\label{ADCeq}
 \sh{A} \sh{C} \subseteq \sD\sD^{\sigma} \cdots \sD^{\sigma^{s-1}}
\eeq
\end{itemize}
 \end{defn}

Given ADC data $\surfD= \surfdata$, we define sheaves $\sh{T}_n$ by
setting $\sh{T}_0 := \struct_X$ and
\[ \sh{T}_n := \sh{A} \sh{D}^{\sigma^s} \cdots \sh{D}^{\sigma^{n-1}} \sh{C}^{\sigma^n} \cdot \Lsh \otimes \Lsh^{\sigma} \otimes \cdots \otimes \Lsh^{\sigma^{n-1}}\]
for $n \geq 1$.
 The sheaves $\sh{T}_n$ satisfy
\[ \sh{T}_n \sh{T}_{m}^{\sigma^n} \subseteq \sh{T}_{n+m},\]
thanks to \eqref{ADCeq}.  Thus we may define
a  $\kk$-algebra 
\[ T (\surfD) :=  \bigoplus_{n \geq 0} H^0(X, \sh{T}_n),\]
where the multiplication is given by
\[\xymatrix{
 H^0(X, \sh{T}_n) \otimes H^0(X, \sh{T}_m) \ar[r]^(.58){1\otimes \sigma^{n}} & H^0(X, \sh{T}_n\otimes \sh{T}_m^{\sigma^n}) \ar[r] & H^0(X, \sh{T}_{n+m}). 
 }\]
We refer to $T$ as an {\em ADC ring}; these rings are our main object of study. 
 
ADC rings generalize classes of  geometric algebras studied previously by the author and by Keeler, Rogalski and Stafford.    In particular, if $\sh{C} = \sh{D} = \struct_X$, then $T(\surfD)$ is the {\em geometric idealizer} $R(X, \Lsh, \sigma, Z)$ studied in \cite{S-idealizer}, where $Z$ is the subscheme defined by $ \sh{A}$.  
If $\sh{C} = \struct_X$ and $\sh{A} = \sD \sD^{\sigma} \cdots \sD^{\sigma^{s-1}}$, then 
$T(\surfD)$ is a na\"ive blowup algebra, as studied in \cite{KRS} and  in \cite{RS-0}.   Of course, if $\sh{ADC} = \struct_X$, then the algebra $T(\surfD)$ is simply the twisted homogeneous coordinate ring $B(X, \Lsh, \sigma)$, as defined in \cite{AV}.  Recall that if $\Lsh$ is appropriately positive -- that is, if $\Lsh$ is {\em $\sigma$-ample} (see Definition~\ref{def-ample}) -- then $B(X, \Lsh, \sigma)$ is noetherian by \cite[Theorem~1.4]{AV} and \cite[Theorem~1.2]{Keeler2000}.

While the definition of the algebras $T(\surfD)$ seems technical, we show in \cite{S-surfclass} that these algebras occur naturally in the classification of noncommutative projective surfaces:  connected $\NN$-graded noetherian domains of GK-dimension 3.  To describe the algebras that occur, we define a geometric condition on ADC data.
 Recall from \cite{S-idealizer} that if $\sigma \in \Aut(X)$, $Z$ is a closed subscheme of $X$, and $A \subset \ZZ$ is infinite, then
$\{ \sigma^n Z \}_{n \in A}$ is   {\em critically transverse} if for all closed subschemes $Y$ of $X$, we have for all but finitely many $n \in A$ that
\[ \shTor_i^X(\struct_{\sigma^n Z}, \struct_Y) = 0\]
for all $i \geq 1$.  (The vanishing of the higher $\shTor$ was called {\em homological transversality} in \cite{S-idealizer}.)  
If $Z$ is 0-dimensional, then $\{\sigma^n(Z)\}_{n \in A}$ is critically transverse if and only if for every $p \in Z$, the set  $\{ \sigma^n(p)\}_{n \in A}$ is  {\em critically dense}.  A subset of $X$ is critically dense if it is infinite, and   any closed subscheme $Y$ contains only finitely many points in the set.

\begin{defn}\label{transverse}
Let $\surfD = \surfdata$ be ADC data.  Let $Z$ be the subscheme of $X$ defined by  $\sh{D}$ and let $\Gamma$ be the subscheme defined by $\sC$.  Using a primary decomposition of $\sA$, write
\beq \label{Omega-def}
\sA = \sI_{\Omega} \cap \sI_{\Lambda},
\eeq
where $\Omega$ is a curve (without embedded components) and $\sI_{\Lambda}$ is maximal with respect to \eqref{Omega-def}.
We say that $\surfD$ is {\em transverse} if the sets
\begin{itemize}
\item $\{\sigma^n Z\}_{n \in \ZZ}$,
\item  $\{ \sigma^n \Omega \}_{n \in \ZZ}$,
\item  $\{ \sigma^n \Lambda\}_{n \geq 0}$, and 
\item $\{ \sigma^n \Gamma \}_{n \leq 0}$ 
\end{itemize}
are critically transverse.  Note that although $\Lambda$ is not uniquely determined by $\sA$, its support is well-defined, and so whether or not $\surfD$ is transverse does not depend on the primary decomposition of $\sA$.
\end{defn}
We note that if $\surfD$ is transverse, then  $\Omega$ is locally principal by Lemma~\ref{lem-CT}, and so the definition of $T(\surfD)$ is left-right symmetric.  That is, there is transverse ADC data $\surfD'$ so that $T(\surfD) \cong T(\surfD')^{op}$.
We also caution the reader that the definitions of ``transverse'' used here and in \cite{S-surfclass} are not precisely equivalent; in particular, in \cite{S-surfclass} we assume in the definition of transversality that $\Lsh$ is $\sigma$-ample.

We now state the main result from \cite{S-surfclass}.

\begin{theorem}\label{ithm-surfclass}
{\em (\cite[Theorem~1.10]{S-surfclass})}
Let $\kk$ be an uncountable algebraically closed field. 
Let $R$ be a connected $\NN$-graded noetherian domain of GK-dimension 3 such that $Q_{\gr}(R) \cong K[z, z^{-1}; \sigma]$ for a field $K$ (necessarily of transcendence degree 2) and automorphism $\sigma$ of $K$.  Then there are transverse ADC data 
\[\surfD = (X, \sL, \sigma, \sA, \sD, \sC, 1),\]
 where $\Lsh$ is $\sigma$-ample, and an integer $k$ so that
\[ R \ver{k} \cong T(\surfD).\]
(Here $R \ver{k}$ denotes the $k$'th Veronese subalgebra of $R$.)
\end{theorem}

The algebras $T(\surfD)$ are therefore clearly of interest.  In particular, it is natural to ask whether transversality is sufficient for $T(\surfD)$ to be noetherian, or if some other, potentially more subtle, property is required.

Our main result is that if $\Lsh$ is $\sigma$-ample, then transversality of the ADC data $\surfD$ is necessary and sufficient for $T(\surfD)$ to be noetherian.  

\begin{theorem}\label{ithm-main}
Let $\surfD = \surfdata$ be  ADC data, where $\Lsh$ is $\sigma$-ample. Then  $T(\surfD)$ is  noetherian if and only if $\surfD$ is transverse.
\end{theorem}

Recall from \cite{S-idealizer} that critical transversality of the set $\{\sigma^n Z \}_{n \geq 0}$ is the controlling property for the geometric idealizer $R(X, \Lsh, \sigma, Z)$ to be noetherian.  Similarly, \cite{KRS} and \cite{RS-0} show that critical density of orbits is the controlling property for na\"ive blowups to be noetherian.  Theorem~\ref{ithm-main} thus generalizes these results.

Let $T$ be an $\NN$-graded $\kk$-algebra.  In noncommutative geometry, one often considers the category
\[ \rqgr T,\]
defined to be the category of graded right $R$-modules modulo torsion.   We analyze this category for  $T=T(\surfD)$, and show:
\begin{theorem}\label{ithm-qgr}
Let $\surfD = \surfdata$ be transverse ADC data where $\Lsh$ is $\sigma$-ample.  Let $T: = T(\surfD)$.  Then:

$(1)$ the functor $\Hom_{\rqgr T}(T, \blank)$ has  finite  cohomological dimension;

$(2)$ $\rqgr T$ depends only on $X, \sigma$, and $\sh{D}$.
\end{theorem}

We also analyze the Artin-Zhang $\chi$ conditions (defined in Section~\ref{CHI}) for the algebra $T(\surfD)$, and determine when $\chi_1$ and $\chi_2$ hold.  In particular, we show in Theorem~\ref{thm-chi2} that  if $\sh{ADC} \neq \struct_X$, then $T(\surfD)$ fails left and right $\chi_2$.  Combined with Theorem~\ref{ithm-surfclass}, this implies that if a birationally commutative surface satisfies left or right $\chi_2$, it is a twisted homogeneous coordinate ring (in sufficiently divisible degree) and satisfies $\chi$.

Let $\surfD = \surfdata$ be ADC data.   If
\[ \sh{A} = (\sD\sD^{\sigma} \cdots \sD^{\sigma^{s-1}}: \sh{C}),\]
we say that the data $\surfD$ is {\em left maximal}. If 
\[ \sh{C} = (\sD\sD^{\sigma} \cdots \sD^{\sigma^{s-1}}: \sh{A}),\]
we say that $\surfD$ is {\em right maximal}.   It is {\em maximal} if is both left and right maximal; that is, if the pair $(\sh{A}, \sh{C})$ is maximal with respect to \eqref{ADCeq}.
ADC rings associated to transverse maximal data are particularly interesting.  
We will see that these rings  have many similar properties to na\"ive blowups at a point, although the algebras are more general.  Further, these algebras give rise to new examples of  maximal orders. 

\begin{theorem}\label{ithm-maxord}
Let $\surfD=(X, \Lsh, \sigma, \sA, \sD, \sC,1)$ be maximal transverse ADC data, where $X$ is a normal surface and $\sL$ is $\sigma$-ample.  Then $T(\surfD)$ is a maximal order.  
\end{theorem}

We summarize the organization of the paper.  In Section~\ref{BACKGROUND}, we recall the definition and basic properties of {\em bimodule algebras}:   roughly speaking, quasicoherent sheaves with a multiplicative structure.  (This is the correct way to work with the sheaf $\bigoplus \sh{T}_n$ defined above.)  In Section~\ref{BASICS}, we give some equivalent formulations of the key condition of transversality of the data $\surfD$, and show that transversality implies that the sheaves $\sh{T}_n$ are ample, in the appropriate sense.  This is a key technical point in proving Theorem~\ref{ithm-main}, which we do in Section~\ref{NOETHERIAN}; we also analyze when the algebras $T(\surfD)$ remain noetherian upon (commutative) base extension.  We study the Artin-Zhang $\chi$ conditions for $T(\surfD)$ in Section~\ref{CHI}.  In Section~\ref{PROJ} we prove Theorem~\ref{ithm-qgr}.  Finally, we prove Theorem~\ref{ithm-maxord} in Section~\ref{MAXORD}.

{\bf Acknowledgements}.  Most of this paper was written while I was a postdoctoral fellow at the University of Washington; I thank the department there, particularly James Zhang and Paul Smith, for their hospitality and for many informative conversations.  Some results  were part of my Ph. D. thesis at the University of Michigan, under the direction of J. T. Stafford.  During the writing of this paper, I was supported by NSF grants   DMS-0555750 and DMS-0802935.  

I am grateful to the referee for many constructive suggestions.

\section{Bimodule algebras}\label{BACKGROUND}

Throughout, we let $\kk$ be a fixed algebraically closed field; all schemes are of finite type over $\kk$.

The subject of this paper is a certain class of graded $\kk$-algebras, defined by geometric data.  As has become standard in the study of subalgebras of twisted homogeneous coordinate rings (see \cite{S-idealizer}, \cite{KRS}, and \cite{RS-0}), one of our main techniques will be to work, not with an algebra, but  with an associated quasi-coherent sheaf on $X$.  This object is known as a {\em bimodule algebra}, and is, roughly speaking, a sheaf with multiplicative structure.   
 In this section, we give the definitions and notation to allow us to work with  bimodule algebras.  Most of the material in this section was developed in \cite{VdB1996} and  \cite{AV}, and we refer the reader there for references.  Our presentation  follows that in \cite{KRS} and \cite{S-idealizer}.

\begin{defn}\label{def-bimod}
Let $X$ be a projective scheme (over $\kk$).  An {\em $\struct_X$-bimodule}  is  a quasicoherent $\struct_{X \times X}$-module $\sh{F}$, such that for every coherent $\sh{F}' \subseteq \sh{F}$,  the projection maps $p_1, p_2: \Supp \sh{F}' \to X$ are both finite morphisms.    The left and right $\struct_X$-module structures associated to an $\struct_X$-bimodule $\sh{F}$ are defined respectively as $(p_1)_* \sh{F}$ and $ (p_2)_* \sh{F}$.  
We make the notational convention that when we refer to an $\struct_X$-bimodule simply as an $\struct_X$-module, we are using the left-handed structure (for example, when we refer to the global sections or higher cohomology
 of an $\struct_X$-bimodule).

There is a tensor product operation on the category of bimodules that has the expected properties; see \cite[Section~2]{VdB1996}.
\end{defn}

All the bimodules that we consider will be constructed from bimodules of the following form:
\begin{defn}\label{def-LR-structure}
Let $X$ be a projective scheme and let $\sigma, \tau \in \Aut(X)$. Let $(\sigma, \tau)$ denote the map
\begin{align*}
X & \to X \times X \\
x & \mapsto (\sigma(x), \tau(x)).
\end{align*}
If  $\sh{F}$ is a quasicoherent sheaf on $X$, we define the $\struct_X$-bimodule ${}_{\sigma} \sh{F}_{\tau}$ to be 
\[ {}_{\sigma} \sh{F}_{\tau} = (\sigma, \tau)_* \sh{F}.\]
If $\sigma = 1$ is the identity, we will often omit it; thus we write $\sh{F}_{\tau}$ for ${}_1 \sh{F}_{\tau}$ and $\sh{F}$ for  the  $\struct_X$-bimodule ${}_1 \sh{F}_1 = \Delta_* \sh{F}$, where $\Delta: X \to X\times X$ is the diagonal.
\end{defn}

The following  lemma  shows how to work with bimodules of the form ${}_{\sigma} \sh{F}_{\tau}$, and, in particular, how to form their tensor product.   
If $\sigma$ is an automorphism of $X$ and $\sh{F}$ is a sheaf on $X$, recall the notation that $\sh{F}^{\sigma} = \sigma^*\sh{F}$.  If $\sL$ is an invertible sheaf on $X$, we define
\[ \sL_n := \sL \otimes \sL^{\sigma} \otimes \cdots \otimes \sL^{\sigma^{n-1}}.\]

\begin{lemma}\label{lem-bimods}
{\em (\cite[Lemma~2.3]{KRS})}
Let $X$ be a projective scheme, let $\sh{F}$, $\sh{G}$ be coherent $\struct_X$-modules, and let $\sigma, \tau \in \Aut(X)$.

$(1)$ ${}_{\tau} \sh{F}_{\sigma} \cong  (\sh{F}^{\tau^{-1}})_{\sigma \tau^{-1}}$.

$(2)$  $\sh{F}_{\sigma} \otimes \sh{G}_{\tau} \cong (\sh{F} \otimes \sh{G}^{\sigma})_{\tau \sigma}$.

$(3)$ In particular, if $\Lsh$ is an invertible sheaf on $X$, then  $\Lsh_{\sigma}^{\otimes n} = (\Lsh_n)_{\sigma^n}$.  \qed
\end{lemma}

\begin{defn}\label{def-BMA}
Let $X$ be a projective scheme and let $\sigma \in \Aut(X)$.  
An {\em $\struct_X$-bimodule algebra}, or simply a {\em bimodule algebra}, $\sh{B}$ is an algebra object in the category of bimodules.  That is, there are a unit map $1: \struct_X \to \sh{B}$ and a product map $\mu: \sh{B} \otimes \sh{B} \to \sh{B}$ that have the usual properties.
\end{defn}

We follow \cite{KRS} and define
\begin{defn}\label{def-gradedBMA}
Let $X$ be a projective scheme and let $\sigma \in \Aut(X)$.  
A bimodule algebra $\sh{B}$ is a {\em graded $(\struct_X, \sigma)$-bimodule algebra} if:

(1) There are coherent sheaves $\sh{B}_n$ on $X$ such that
\[ \sh{B} = \bigoplus_{n \in \ZZ} {}_1(\sh{B}_n)_{\sigma^n};\]

(2)  $\sh{B}_0 = \struct_X $;

(3) the multiplication map $\mu$ is given by $\struct_X$-module maps
$\sh{B}_n \otimes \sh{B}_m^{\sigma^n} \to \sh{B}_{n+m}$, satisfying the obvious associativity conditions.  
\end{defn}

\begin{defn}\label{def-module}
Let $X$ be a projective  scheme and let $\sigma \in \Aut(X)$.  
Let $\sh{R} = \bigoplus_{n\in\ZZ} (\sh{R}_n)_{\sigma^n}$ be a graded $(\struct_X, \sigma)$-bimodule algebra.  A {\em right $\sh{R}$-module} $\sh{M}$ is a quasicoherent $\struct_X$-module $\sh{M}$ together with a right $\struct_X$-module map $\mu: \sh{M} \otimes \sh{R} \to \sh{M}$ satisfying the usual axioms.  We say that $\sh{M}$ is {\em graded} if there is a direct sum decomposition 
\[\sh{M} = \bigoplus_{n \in \ZZ}  (\sh{M}_n)_{\sigma^n}\]
 with multiplication giving a family of $\struct_X$-module maps $\sh{M}_n \otimes \sh{R}_m^{\sigma^n} \to \sh{M}_{n+m}$, obeying the appropriate axioms.
 
  We say that $\sh{M}$ is {\em coherent} if there are a coherent $\struct_X$-module $\sh{M}'$ and a surjective map $\sh{M}' \otimes \sh{R} \to \sh{M}$ of ungraded $\struct_X$-modules.  
We make similar definitions for left $\sh{R}$-modules.  The bimodule algebra $\sh{R}$ is {\em right (left) noetherian}
if every right (left) ideal of $\sh{R}$ is coherent.  A graded $(\struct_X, \sigma)$-bimodule algebra is right (left) noetherian if and only if every graded right (left) ideal is coherent.
\end{defn}

We recall here some standard notation for module categories over rings and bimodule algebras.  Let $R$ be an $\NN$-graded $\kk$-algebra.  We define $\rGr R$ to be the category of $\ZZ$-graded right $R$-modules; morphisms in $\rGr R$ preserve degree.   Let $\rTors R$ be the full subcategory of modules that are direct limits of right bounded modules.  This is a Serre subcategory of $\rGr R$, so we may form the {\em quotient category}
\[ \rQgr R := \rGr R / \rTors R.\]
(We refer the reader to \cite{Gabriel1962} as a reference for the category theory used here; note that the convention there, which we follow, is that objects of  $\rGr R$ are also objects of $\rQgr R$.)  There is a canonical  quotient functor from $\rGr R$  to $ \rQgr R$.  

We make similar definitions on the left.   Further, throughout this paper, we adopt the convention that if Xyz is a category, then xyz is the full subcategory of noetherian objects.  Thus we have $\rgr R$ and $\rqgr R$, $R \lqgr$, etc.  If $X$ is  a scheme, we will denote the category of quasicoherent (respectively coherent) sheaves on $X$ by $\struct_X \lMod$ (respectively $\struct_X \lmod$).  

Given a module $M \in \rgr R$, we define $M[n] = \bigoplus_{i \in \ZZ} M[n]_i$, where
\[ M[n]_i = M_{n+i}.\]
If $M, N \in \rgr R$, let 
\[ \underline{\Hom}_{\rgr R} (M, N) = \bigoplus_{n \in \ZZ} \Hom_{\rgr R}(M, N [n]).\]
Similarly, if $\sh{M}, \sh{N} \in \rqgr R$, we define
\[ \underline{\Hom}_{\rqgr R} (\sh{M}, \sh{N}) = \bigoplus_{n \in \ZZ} \Hom_{\rqgr R}(\sh{M}, \sh{N} [n]).\]
The $\underline{\Hom}$ functors have derived functors  $\underline{\Ext}_{\rgr R}$ and $\underline{\Ext}_{\rqgr R}$.  

For a graded $(\struct_X, \sigma)$-bimodule algebra $\sh{R}$, we likewise define $\rGr \sh{R}$ and $\rgr \sh{R}$.  The full subcategory  $\rTors \sh{R}$ of $\rGr \sh{R}$ consists of direct limits of modules that are  coherent as $\struct_X$-modules, and we similarly define
\[ \rQgr \sh{R} := \rGr \sh{R}/\rTors \sh{R}.\]
We define $\rqgr \sh{R}$ in the obvious way.  

If $R$ is a graded domain, then a graded right $R$-module $M$ is {\em Goldie torsion} if any homogeneous $m \in M$ is annihilated by some nonzero homogeneous $r\in R$; equivalently, $M$ is a direct limit of sums of modules of the form $(R/I)[n]$ for some graded right ideal $I$ of $R$.  If $X$ is a projective variety and $\sh{R} $ is a graded $(\struct_X, \sigma)$-bimodule algebra, we say that a graded right  $\sh{R}$-module $\sh{M}$ is Goldie torsion if $\sh{M}$ is a direct limit of sums of modules of the form $(\sh{R}/\sh{I})[n]$ for a graded right ideal $\sh{I}$ of $\sh{R}$.  We denote the full subcategory of $\rgr R$ (respectively, $\rgr \sh{R}$) consisting of Goldie torsion modules by $\GT(\rgr R)$ (respectively, $\GT(\rgr \sh{R})$).

If $\sh{R}$ is an $\struct_X$-bimodule algebra, its global sections $H^0(X, \sh{R})$ inherit a $\kk$-algebra structure.  We call $H^0(X, \sR)$ the {\em section algebra} of $\sR$.  If $\sR = \bigoplus (\sR_n)_{\sigma^n}$ is a graded $(\sO, \sigma)$-bimodule algebra, then multiplication on $H^0 (X, \sh{R})$ is induced from  the maps
\[ \xymatrix{
H^0(X, \sh{R}_n) \otimes H^0(X, \sh{R}_m) \ar[r]^{1\otimes \sigma^n} 
&  H^0(X, \sh{R}_n) \otimes H^0(X, \sh{R}_m^{\sigma^n}) \ar[r]^<<<<<{\mu} 
& H^0(X, \sh{R}_{n+m}).
}\]

  If $\sh{M}$ is a graded right $\sh{R}$-module, then 
\[ H^0(X, \sh{M}) = \bigoplus_{n \in \ZZ} H^0(X, \sh{M}_n)\]  
 is a right $H^0(X, \sh{R})$-module in the obvious way; thus $H^0(X, \blank)$ is a functor from $\rGr \sh{R}$ to $\rGr H^0(X, \sh{R})$.  

 If  $R = H^0(X, \sh{R})$ and $M$ is  a graded right $R$-module, define 
$M \otimes_R \sh{R}$ to be the sheaf associated to the presheaf $V \mapsto M \otimes_R \sh{R}(V)$.   This is  a graded right $\sh{R}$-module, and the functor $\blank \otimes_R \sh{R}: \rGr R \to \rGr \sh{R}$ is a right adjoint to $H^0(X, \blank)$.

\begin{defn}\label{def-ample}
Let $X$ be a projective scheme, let $\sigma \in \Aut(X)$, and let  $\{\sh{R}_n\}_{n \in \NN}$ be a sequence of coherent sheaves  on $X$.  The sequence of bimodules $\{(\sh{R}_n)_{\sigma^n}\}_{n \in \NN}$ is {\em right ample} if for any coherent $\struct_X$-module $\sh{F}$, the following properties hold:
 
 (i)   $\sh{F} \otimes \sh{R}_n$ is globally generated for $n \gg 0$;
 
 (ii) $H^q(X, \sh{F} \otimes \sh{R}_n) = 0$  for $n \gg 0$ and $q \geq 1$.  
 
\noindent The sequence $\{ (\sh{R}_n)_{\sigma^n}\}_{n \in \NN}$  is {\em left ample}  if for any coherent $\struct_X$-module $\sh{F}$, the following properties hold:
 
 (i)  $\sh{R}_n \otimes \sh{F}^{\sigma^n}$ is globally generated for $n \gg 0$;
 
 (ii) $H^q(X, \sh{R}_n \otimes \sh{F}^{\sigma^n}) = 0$  for $n \gg 0$ and  $q \geq 1$.

We say that an invertible  sheaf $\Lsh$ is {\em $\sigma$-ample}
 if the $\struct_X$-bimodules 
\[\{(\Lsh_n)_{\sigma^n}\}_{n \in \NN} = \{ \Lsh_\sigma^{\otimes n}\}_{n \in \NN}\]
 form a right ample sequence.  By  \cite[Theorem~1.2]{Keeler2000}, this is true if and only if the $\struct_X$-bimodules $\{(\Lsh_n)_{\sigma^n}\}_{n \in \NN}$ form a left ample sequence.
\end{defn}

The following result is a special case of a result due to Van den Bergh \cite[Theorem~5.2]{VdB1996}, although we follow the presentation of  \cite[Theorem~2.12]{KRS}:

\begin{theorem}\label{thm-VdBSerre}
{\em (Van den Bergh)}
Let $X$ be a projective scheme and let $\sigma$ be an automorphism of $X$.  Let $\sh{R} = \bigoplus (\sh{R}_n)_{\sigma^n}$ be a right noetherian graded $(\struct_X, \sigma)$-bimodule algebra, such that the bimodules  $\{(\sh{R}_n)_{\sigma^n}\}$ form a right ample sequence.   Then $R = H^0(X, \sh{R})$ is also right noetherian, and the functors $H^0(X, \blank)$ and $\blank \otimes_R \sh{R}$ induce an equivalence of categories
\[ \rqgr \sh{R} \simeq \rqgr R.\] 
\qed
\end{theorem}

The fundamental example of a bimodule algebra is the following.  
Let $X$ be a projective scheme, let $\sigma \in \Aut(X)$, and let $\Lsh$ be an invertible sheaf on $X$.  We define the {\em twisted bimodule algebra of $\Lsh$} to be
\[ \sh{B} = \sh{B}(X, \Lsh, \sigma) = \bigoplus_{n \geq 0} ( \Lsh_n)_{\sigma^n}.\]
Then $\sh{B}$ is an $(\struct_X, \sigma)$-graded bimodule algebra.  
Taking global sections of $\sh{B}(X, \Lsh, \sigma)$ gives the twisted homogeneous coordinate ring $B(X, \Lsh, \sigma)$.

Throughout this paper, we will consider sub-bimodule algebras of 
the twisted bimodule algebra $\sh{B} =\sh{B}(X, \Lsh, \sigma)$.  We note here that  the invertible sheaf $\Lsh$ makes only a formal difference.  

\begin{lemma}\label{lem-catequiv}
{\em (\cite[Lemma~2.12]{S-idealizer})}
Let $X$ be a projective scheme with automorphism $\sigma$, and let $\Lsh$ be an invertible sheaf on $X$.  Let 
\[ \sh{R} = \bigoplus_{n \geq 0} (\sh{R}_n)_{\sigma^n}\]
 be a graded $(\struct_X,\sigma)$-sub-bimodule algebra of the twisted bimodule algebra $\sh{B}(X, \Lsh, \sigma)$.  Let $\sh{S}_n := \sh{R}_n \otimes  \Lsh_n^{-1}$ for $n \geq 0$.   

  Let $\sh{S}$ be the graded $(\struct_X, \sigma)$-bimodule algebra defined by 
 \[ \sh{S} := \bigoplus_{n \geq 0} (\sh{S}_n)_{\sigma^n} .\]
 Then the categories $\rgr \sh{R}$ and $\rgr \sh{S}$ are equivalent, and the categories $\sh{S} \lgr$ and $\sh{R} \lgr $ are equivalent. 
\qed
\end{lemma}

 To end the section, we record the effect of shifting degrees on a graded $(\struct_X, \sigma)$-bimodule algebra.

\begin{lemma}\label{lem-shift}
Let $X$ be a projective scheme, let $\sigma \in \Aut(X)$, and let $\Lsh$ be an invertible sheaf on $X$.  Let $\sh{R} = \bigoplus_n (\sh{R}_n)_{\sigma^n}$ be a graded sub-$(\struct_X, \sigma)$-bimodule algebra of the twisted bimodule algebra $\sh{B}(X, \Lsh, \sigma)$, and let $\sh{N}$ be a graded right $\sh{R}$-module.  We write
\[ \sh{N}_n = \sh{F}_n \otimes \Lsh_n,\]
where $\sh{F}_n$ is a quasicoherent sheaf on $X$ with trivial bimodule structure.  
 If $m \in \ZZ$, then $\sh{N} [m] \cong \bigoplus \sh{G}_n \otimes \Lsh_{\sigma}^{\otimes n}$, where:

$\sh{G}_n = (\sh{F}_{n+m} \otimes \Lsh_m)^{\sigma^{-m}} $ if $m > 0$ (with the trivial bimodule structure), and

$\sh{G}_n = (\sh{F}_{n+m})^{\sigma^{-m}} \otimes \Lsh_{-m}^{-1}$ if $m < 0$.
\end{lemma}
\begin{proof}
This follows exactly as in the proof of \cite[Lemma~5.5]{KRS}.
\end{proof}

\section{Basic properties of geometric algebras}\label{BASICS}

Let $\surfD = \surfdata$ be ADC data, as defined in the Introduction.  Recall that we define  sheaves $\sh{T}_n$ by setting  $\sh{T}_0 := \struct _X$ and 
\[\sh{T}_n:=  \sh{A} \sh{D}^{\sigma^s}\cdots \sh{D}^{\sigma^{n-1}}\sh{C}^{\sigma^n}\Lsh_n\]
for $n \geq 1$.  (Note that for $n \leq s$ we have $\sT_n = \sA \sC^{\sigma^n} \sL_n$.)
We then define a bimodule algebra 
\[ \sh{T}(\surfD):= \bigoplus_n (\sh{T}_n)_{\sigma^n}.\]
The ring $T(\surfD)$ is thus the section ring of $\sh{T}(\surfD)$.   

Given ADC data $\surfD = \surfdata$, let $Z$ be the subscheme of $X$ defined by $\sD$, and let $\Gamma$ be the subscheme defined by $\sC$.  Let 
\[ \sA = \sP_1 \cap \cdots \cap \sP_k \cap \sQ_1 \cap \cdots \sQ_{\ell}\]
be a minimal primary decomposition of $\sA$, where the $\sP_i$ have height 1 associated primes and the $\sQ_j$ have maximal associated primes.  Let $\Omega$ be the curve defined by 
\[ \sI_{\Omega} = \sP_1 \cap \cdots \cap \sP_k\]
and let $\Lambda$ be the 0-dimensional subcheme defined by 
\[ \sI_{\Lambda} = \sQ_1 \cap \cdots \cap \sQ_{\ell}.\]
  We then have
\beq \label{Omega-def2}
\sA = \sI_{\Omega} \cap \sI_{\Lambda},
\eeq  
We call the tuple
\[ \dot{\surfD} := \geomdata\]
the {\em geometric data associated to $\surfD$}.
 
Recall from the introduction that $\surfD$ is {\em transverse} if the sets
 $\{ \sigma^n \Omega\}_{n \in \ZZ}$, 
  $\{\sigma^n Z\}_{n \in \ZZ}$,
 $\{ \sigma^n \Lambda\}_{n \geq 0}$, and 
 $\{ \sigma^n \Gamma \}_{n \leq 0}$ 
are critically transverse.  
Note that $\Lambda$ is not defined uniquely, but that the support of $\Lambda$ is well-defined; in particular, the transversality of $\surfD$ does not depend on the choice of $\Lambda$.  If  $\surfD$ is transverse, then both $\Lambda$ and $\Gamma$ are, in particular, supported on infinite orbits.  Note also that in this case, if we are willing to replace $T(\surfD)$ by a Veronese subring, we may always assume that $s=1$.

Ultimately, we will show that transversality of $\surfD$ implies that $T(\surfD)$ is noetherian.  In this section, we analyze the definition of transversality and give simpler equivalent formulations.  We then study when the bimodules  $\{(\sh{T}_n)_{\sigma^n} \}$ form an ample sequence.

  We note that on a surface, we may reframe the condition for critical transversality of the $\sigma$-orbit of a curve.

\begin{lemma}\label{lem-CT}
Let $X$ be a projective scheme, and let $\sigma \in \Aut(X)$.  Let $\Omega \subseteq X$ be a closed subscheme of pure codimension 1.  The following are equivalent:

$(1)$ $\{\sigma^n \Omega\}_{n \in \ZZ}$ is critically transverse;

$(2)$ $\Omega$ contains no reduced and irreducible subschemes that are of finite order under $\sigma$,  meets orbits only finitely often, and is locally principal.
\end{lemma}
\begin{proof}
Suppose that $\Omega$ is locally principal and that $W$ is a reduced and irreducible proper subscheme of $X$.  We claim that $W \subseteq \Omega$ if and only if $\shTor^X_1(\struct_{\Omega}, \struct_{W}) \neq 0$.

To prove the claim it suffices to work locally.  So let $x \in \Omega \cap W$.  Let $A := \struct_{X,x}$, let $P$ be the prime ideal defining $W$ in $A$, and let $a \in A$ be the local equation of $\Omega$. By \cite[Exercise~3.1.3]{Weibel}, we may identify 
\beq\label{penguin}
 \Tor^A_1(A/aA, A/P) \cong (aA \cap P)/aP.
\eeq
 Let $J := (P:  aA)$, so $aA \cap P = a J$.  Note $J \supseteq P$. Then $J \supsetneqq P$ if and only if $aA \cap P \neq aP$, which by \eqref{penguin} happens exactly when $\Tor^A_1(A/aA, A/P) \neq 0$.

Suppose that $a \not\in P$.  Since $aJ \subseteq P$ and $P$ is prime, we must have $J=P$.    On the other hand, if $a \in P$, then $J = A \supsetneqq P$.  Thus $\Tor^A_1(A/aA, A/P) \neq 0 $ if and only if $a \in P$, as claimed.

$(2) \Rightarrow (1)$.  
Assume that (2) holds.  
By \cite[Lemma~5.7]{S-idealizer}, to show that $\{\sigma^n \Omega\}$ is critically transverse, it suffices to prove that for all reduced and irreducible subschemes $Y$ of $X$, the set
\[ \{ n \in \ZZ \st \shTor^X_1(\struct_{\sigma^n \Omega}, \struct_{Y}) \neq 0 \}\]
is finite.  This follows directly from the claim above, since by assumption the set
\[ \{ n \in \ZZ \st Y \subseteq  \sigma^n \Omega\} \] 
is finite for any reduced and irreducible $Y$.

$(1) \Rightarrow (2)$.
Suppose that $\{\sigma^n \Omega\} $ is critically transverse.  By \cite[Lemma~7.7]{S-idealizer}, $\Omega$ is locally principal.  The rest of $(2)$ is immediate  from the claim.
\end{proof}

In characteristic 0, the conditions for transversality of ADC data simplify even further.
\begin{proposition}\label{prop-transverse-char0}
Let $\kk$ have characteristic 0, and let $\surfD = \surfdata$ be ADC data, with associated geometric data $\dot{\surfD} = \geomdata$.  Then the following are equivalent:

$(1)$  $\surfD$ is transverse;

$(2)$ $\Omega$ is locally principal and contains no points or components of finite order under $\sigma$, and all points in $\Lambda \cup Z \cup \Gamma$ have dense orbits.
\end{proposition}
\begin{proof}
$(1) \Rightarrow (2)$ follows from Lemma~\ref{lem-CT}.

$(2) \Rightarrow (1)$.  
That $\{\sigma^n \Lambda\}_{n \geq 0}$, $\{\sigma^n \Gamma\}_{n \leq 0}$, and $\{\sigma^n Z\}_{n \in \ZZ}$ are critically transverse follows from \cite[Theorem~5.1]{BGT}.  Suppose that $\Omega \cap \{ \sigma^n( x)\}$ is infinite for some $x \in X$.  By \cite[Theorem~5.1]{BGT}, $\{\sigma^n (x)\}$ is not Zariski-dense in $X$.  Thus 
\[ \bbar{\{ \sigma^n (x)\}} = C_1 \cup C_2 \cup \cdots \cup C_k \]
consists of finitely many irreducible curves, which are trivially of finite order under $\sigma$.  Some $C_i$ therefore meets $\Omega$ infinitely often and is thus contained in $\Omega$, a contradiction.  Lemma~\ref{lem-CT} now implies that $\{ \sigma^n \Omega\}$ is critically transverse and $\surfD$ is transverse.
\end{proof}

We note that, while the definition of the ring $T(\surfD)$ has left-right asymmetry, if $\surfD$ is transverse then in fact the definition is symmetric.  Since $\Omega$ is locally principal, we may let
\[ \Lsh' := \Lsh(-\Omega+\sigma^{-1}(\Omega)).\]
Then $\Lsh'$ is also $\sigma$-ample, and therefore $\sigma^{-1}$-ample by \cite[Theorem~1.2]{Keeler2000}.
Define
\[ \sA' := (\sI_{\Omega})^{-1} \sA.\]
Let 
\[ \sC':= \sI_{\Omega} \sC\]
and let $\sD' := \sD^{\sigma^{s-1}}$.
  Then
\[ \sh{T}_n = \sh{I}_{\Omega} \sh{A}'\sh{D}^{\sigma^s} \cdots \sh{D}^{\sigma^{n-1}} \sh{C}^{\sigma^n} \Lsh_n 
\cong \sh{A}(\sh{D}')^{\sigma} \cdots (\sh{D}')^{\sigma^{n-s}} (\sC')^{\sigma^n} (\Lsh')_n.\]
Let
\[\surfD':= (X, \Lsh', \sigma^{-1}, \sh{C}', \sh{D}', \sh{A}',s).\]
Then $\surfD'$ is transverse ADC data, and $T(\surfD) \cong T(\surfD')^{op}$.

Let $\surfD = \surfdata$ be transverse  ADC data and let  
$\sh{T} := \sh{T} (\surfD)$.
   To end the  section, we show  that the sequence of bimodules $\{ (\sh{T}_n)_{\sigma^n}\}$ is left and right ample, in the sense of Definition~\ref{def-ample}.    
We will use a lemma of Rogalski and Stafford  that relates the ampleness of a sequence of bimodules of the form $\{ (\sh{R}_n)_{\sigma^n}\}$ to the Castelnuovo-Mumford regularity of the sheaves $\sh{R}_n$.  (We refer the reader to \cite{RS-0} for the definition of Castelnuovo-Mumford regularity.)  

We recall two results that we will use.  
\begin{lemma}\label{lem-reg}
{\em (\cite[Corollary~3.14]{RS-0}) }
Let $X$ be a projective scheme with very ample invertible sheaf $\sh{N}$.  Let $\sh{F}_n$ be a sequence of coherent sheaves on $X$ such that for each $n$, the closed set where $\sh{F}_n$ is not locally free has dimension at most 2.  Then $\{(\sh{F}_n)_{\sigma^n}\}$ is a right ample sequence if and only if 
\[\lim_{n \to \infty} \reg_{\sh{N}}\sh{F}_n = - \infty,\]
 and $\{ ( \sh{F}_n)_{\sigma^n} \}$ is a left ample sequence if and only if 
 \[\lim_{n \to \infty} \reg_{\sh{N}^{\sigma^n}} \sh{F}_n = - \infty.\]
\end{lemma}
\begin{proof}
The right ampleness statement is a restatement of \cite[Corollary~3.14]{RS-0}.  The left ampleness statement follows by symmetry.
\end{proof}

\begin{lemma}\label{lem-Dennis}
{\em (\cite[Proposition~2.8]{Keeler2006})}
Let $X$ be a projective scheme with very ample invertible sheaf $\sh{N}$.  Then there is a constant $C$, depending only on $X$ and $\sh{N}$, so that for any pair $\sh{F}, \sh{G}$ of coherent sheaves such that the dimension of the closed set where both $\sh{F}$ and $\sh{G}$ are not locally free is less than or equal to 2, we have that
\[ \reg_{\sh{N}} \sh{F} \otimes \sh{G} \leq \reg_{\sh{N}} \sh{F} + \reg_{\sh{N}} \sh{G} + C.\]
\qed
\end{lemma}

We will also frequently use the following easy observation about cohomology vanishing.
\begin{lemma}\label{lem-FOOBAR}
Let $X$ be a projective scheme and suppose that 
\[ 0 \to \sh{K} \to \sh{M} \stackrel{\theta}{\to} \sh{N} \to \sh{K}' \to 0\]
is an exact sequence of coherent sheaves on $X$, where $\sh{K}$ and $\sh{K'}$ are supported on subschemes of dimension 0.  Further suppose that $H^i(X, \sh{M}) = 0$ for all $i \geq 1$.  Then $H^i(X, \sh{N}) = 0 $ for all $i \geq 1$.
\end{lemma}
\begin{proof}
Note that $H^i(X, \sh{K}) = H^i(X, \sh{K}')  = 0$ for all $i \geq 1$.  Let $\sh{M}' := \im \theta$.  From the long exact cohomology sequence, we deduce that $H^i(X, \sh{M}')  =0$ for all $i \geq 1$.  This implies that $H^i(X, \sh{N}) = 0$ for  all $i\geq 1$.
\end{proof}

\begin{corollary}\label{cor-FOOBAR}
Let $X$ be a projective scheme and let $\sigma \in \Aut(X)$.  Let $\sh{M}_n$, $\sh{N}_n$ be a sequence of coherent sheaves on $X$, and suppose there are exact sequences
\[ 0 \to \sh{K}_n \to \sh{M}_n \to \sh{N}_n \to \sh{K}'_n \to 0,\]
where $\sh{K}_n$ and $\sh{K}'_n$ are supported on sets of dimension 0.  Assume that
$\{(\sh{M}_n)_{\sigma^n}\}$ is left (right) ample.  Then $\{(\sh{N}_n)_{\sigma^n}\}$ is left (right) ample.
\end{corollary}
\begin{proof}
We prove the right ampleness statement.  By Lemma~\ref{lem-FOOBAR}, we have for any coherent $\sh{F}$ that $H^i(X, \sh{F}\otimes \sh{N}_n) = 0$ for $i > 0$ and $n \gg 0$.  It follows as in the proof of \cite[Lemma~4.2]{KRS} that $\sh{F}\otimes \sh{N}_n$ is globally generated for $n \gg 0$.
\end{proof}

We will show that the sequence of bimodules $\{ (\sh{T}_n)_{\sigma^n} \}$ is left and right ample under slightly less restrictive assumptions on the defining data than transversality.   

 
\begin{lemma}\label{lem-ample1}

$(1)$  Let $X$ be a projective surface, let $\sigma \in \Aut(X)$, and let $\Lsh$ be a $\sigma$-ample invertible sheaf on $X$.  Let $\Omega$ be a curve on $X$ so that $\{ \sigma^n{\Omega}\}$ is critically transverse.  Let $\sh{E}$ be an ideal sheaf on $X$ that defines a 0-dimensional subscheme supported on dense orbits.  Then the sequence of bimodules 
\[ \{ \bigl( \sh{I}_{\Omega}  \sh{E} \sh{E}^{\sigma} \cdots \sh{E}^{\sigma^{n-1}} \Lsh_n \bigr)_{\sigma^n} \} \]
is left and right ample.

$(2)$  Let $\surfD = \surfdata$ be ADC data, and let $\dot{\surfD} = \geomdata$ be the associated geometric data.  Suppose  $\Lsh$ is $\sigma$-ample, $\{\sigma^n \Omega\}$ is critically transverse, and that all points in $\Lambda \cup Z \cup \Gamma$ lie on dense $\sigma$-orbits.  
Let $\sh{T} := \sh{T}(\surfD)$.  
   Then the sequence of bimodules $ \{(\sh{T}_n)_{\sigma^n}\}$ is left and right ample.
  \end{lemma}
  \begin{proof}
  (1) For all $n \geq 1$, let 
  \[ \sh{J}_n :=  \sh{I}_{\Omega}  \sh{E}  \sh{E}^{\sigma} \cdots \sh{E}^{\sigma^{n-1}} .\]
  We will show that the sequence $\{ (\sh{J}_n \Lsh_n)_{\sigma^n}\}$ is left and right ample.  
   
  We first assume in addition that $\Lsh$ is ample.  By \cite[Theorem~1.7]{AV}, $\Lsh$ is then also $\sigma^2$-ample.  Note that all points in the cosupport of $\sh{E} \sh{E}^{\sigma}$ have dense $\sigma^2$-orbits.  
   Let  
\begin{multline*}
\sh{F}_n := (\sh{E} \sh{E}^{\sigma})(\sh{E} \sh{E}^{\sigma})^{\sigma^2} \cdots (\sh{E} \sh{E}^{\sigma})^{\sigma^{2n-2}} \Lsh \otimes \Lsh^{\sigma^2} \otimes \cdots \otimes \Lsh^{\sigma^{2n-2}} \\
= \sh{E} \sh{E}^{\sigma} \cdots  \sh{E}^{\sigma^{2n-1}} \Lsh \otimes \Lsh^{\sigma^2} \otimes \cdots \otimes \Lsh^{\sigma^{2n-2}}.
\end{multline*}
 By \cite[Theorem~3.1]{RS-0}, the sequences $\{ ( \sh{F}_n)_{\sigma^{2n}} \}$ and $\{ (\sh{F}_{n+1})_{\sigma^{2n+1} }\}$
 are  left and right ample.   
 
 Now let 
 \[ \sh{G}_n := \sh{I}_{\Omega} \Lsh^{\sigma} \Lsh^{\sigma^3} \cdots \Lsh^{\sigma^{2n-1}}.\]
 The sequences
 $\{ (\sh{G}_n)_{\sigma^{2n}}\}$ and $\{ (\sh{G}_n)_{\sigma^{2n+1}}\}$ are left and right ample by \cite[Lemma~6.1 and Proposition~6.2]{S-idealizer}.
 By Lemma~\ref{lem-Dennis} and Lemma~\ref{lem-reg}, the sequences
 \[ \{ (\sh{F}_n \otimes \sh{G}_n)_{\sigma^{2n}} \}\]
 and
 \[ \{ (\sh{F}_{n+1} \otimes \sh{G}_{n})_{\sigma^{2n+1}} \} \]
 are left and right ample.  
 
  For any $n \geq 0$, there is an exact sequence
 \[ 0 \to \sh{H}_n \to  \sh{F}_n \otimes  \sh{G}_n  \to
 \sh{J}_{2n} \Lsh_{2n}  \to 0\]
 where $\sh{H}_n$ is supported on a dimension 0 subscheme of $X$.  By Corollary~\ref{cor-FOOBAR},  $ \{ ( \sh{J}_{2n} \Lsh_{2n})_{\sigma^{2n}}\} $ is  a left and right ample sequence.  Likewise, from the maps
 \[  \sh{F}_{n+1} \otimes \sh{G}_n \to \sh{J}_{2n+1} \Lsh_{2n+1} \]
 we obtain that $\{ (\sh{J}_{2n+1} \Lsh_{2n+1})_{\sigma^{2n+1}}\}$ is left and right ample.  Thus 
 \[ \{ (\sh{J}_n  \Lsh_n )_{\sigma^n}  \} \]
 is left and right ample.
 
  Now consider the general case.  By \cite[Theorem~1.7]{AV}, there is some $k \geq 1$ so that $\Lsh_k$ is ample.  Let $\sh{E}' := \sh{E} \sh{E}^{\sigma} \cdots \sh{E}^{\sigma^{k-1}}$.  We have seen that the sequence of bimodules 
  \[ \{ \bigl(  \sh{I}_{\Omega}  \sh{E}' (\sh{E}')^{\sigma^k} \cdots (\sh{E}')^{\sigma^{k(n-1)}} \Lsh_{kn}\bigr)_{\sigma^{kn}}\} = \{ (\sh{J}_{kn} \Lsh_{kn})_{\sigma^{kn}} \}\]
  is left and right ample.  
    Lemma~\ref{lem-reg} implies that for any $0 \leq i \leq k-1$, the sequence
    \[ \{ ( \sh{J}_{kn} \Lsh_{kn-i})_{\sigma^{kn-i}} \}\]
is left and  right ample.  

Fix $0 \leq i \leq k-1$.  We have $\sh{J}_{kn} \subseteq \sh{J}_{kn-i}$ for all $n \geq 1$, and the factor is supported on  a set of dimension 0.   Thus by Corollary~\ref{cor-FOOBAR} the sequence
\[ \{ ( \sh{J}_{kn-i} \Lsh_{kn-i})_{\sigma^{kn-i}} \}_{n \geq 0} \]
is left and right ample for all $0 \leq i \leq k-1$.  Thus 
\[ \{ ( \sh{J}_n \Lsh_n)_{\sigma^n} \} \]
is a left and right ample sequence, as claimed.

(2)  Let $\sh{E} :=  \sI_{\Lambda} \sh{D}  \sh{C}^{\sigma}$, so for $n \geq 1$
\[ \sI_{\Omega} \sh{E} \sh{E}^{\sigma} \cdots \sh{E}^{\sigma^{n-1}} \subseteq  \sI_{\Omega} \sI_{\Lambda} \sh{D}^{\sigma^s} \cdots \sh{D}^{\sigma^{n-1}} \sh{C}^{\sigma^n} \subseteq  \sh{A} \sh{D}^{\sigma^s} \cdots \sh{D}^{\sigma^{n-1}} \sh{C}^{\sigma^n}.\]
The cokernel of this inclusion is supported on a set of dimension 0.  By (1) the sequence
\[ \{  (\sh{I}_{\Omega}  \sh{E} \sh{E}^{\sigma} \cdots \sh{E}^{\sigma^{n-1}} \Lsh_n )_{\sigma^n} \} \]
is left and right ample.  As above,
\[ \{ (\sh{T}_n)_{\sigma^n} \} \]
is left and right ample.
   \end{proof}

\section{Noetherian rings and bimodule algebras}\label{NOETHERIAN}

We will now prove that if the  ADC data $\surfD = \surfdata$ is transverse and $\Lsh$ is $\sigma$-ample, then both the bimodule algebra  $\sh{T}(\surfD)$ and the $\kk$-algebra $T (\surfD)$  are left and right noetherian. 
We also prove that the converse holds when $\Lsh$ is $\sigma$-ample, and  analyze when $T(\surfD)$ is strongly noetherian.

These proofs are carried out in several steps.  We first analyze the case of maximal transverse ADC data.    
To show that ADC bimodule algebras of maximal transverse data are noetherian, we explicitly construct generators for graded right and left ideals.  Note that if $\surfD = \surfdata$ is maximal ADC data, then $\sD\cdots \sD^{\sigma^{s-1}} \subseteq \sA \cap \sC$; in particular, the cosupport of $\sA$ is 0-dimensional.

\begin{proposition}\label{prop-circ2}
Suppose that the tuple $\ADCD = \surfdata$ is maximal  transverse ADC data, and let $\sh{S} := \sh{T}(\ADCD)$.    

$(1)$
Let $\sh{J} = \bigoplus (\sh{J}_n)_{\sigma^n}$ be a graded right ideal of $\sh{S}$.  Then there are an integer $m \geq s$ and an ideal sheaf $\sh{J}' \subseteq \sh{A} \sh{D}^{\sigma^s} \cdots \sh{D}^{\sigma^{m-1}}$ on $X$   so that $\sJ'$ and $\sD^{\sigma^n}$ are comaximal for $n \geq m$, and for $n \geq m$, 
\[ \sh{J}_n = (\sh{J}' \sh{D}^{\sigma^{m}} \cdots \sh{D}^{\sigma^{n-1}} \sh{C}^{\sigma^n}) \Lsh_n
= \sh{J}'\Lsh_n \cap \sh{S}_n
.\]

$(2)$  Let $\sh{K}$ be a graded left ideal of $\sh{S}$.  Then there are an integer $m' \geq s$ and an ideal sheaf $\sh{K}' \subseteq \sh{D}^{\sigma^{-m'+1}}\cdots\sh{D}^{\sigma^{-1}}\sh{C}$ on $X$   so that $\sK'$ and $\sD^{\sigma^n}$ are comximal for $n \leq -m'$, and for $n \geq m' +s$, 
\[ \sh{K}_n = (\sh{A} \sh{D}^{\sigma^s} \cdots \sh{D}^{\sigma^{n-m'}} (\sh{K}')^{\sigma^n}) \Lsh_n
= \sh{S}_n \cap (\sh{K}')^{\sigma^n} \Lsh_n.\]

$(3)$  Let $\sh{H}$ be a graded ideal of $\sh{S}$.  Then there are a $\sigma$-invariant ideal sheaf $\sh{H}'$ on $X$ and an integer $m'' \geq s$ so that for $n \geq m''$, 
\[ \sh{H}_n = \sh{H}' \sh{S}_n = \sh{H}' \Lsh_n \cap \sh{S}_n.\]
\end{proposition}

\begin{proof}  
   By Lemma~\ref{lem-catequiv}, without loss of generality we may assume that $\Lsh = \struct_X$.
Let $\dot{\surfD}=\geomdata$ be the geometric data associated to $\surfD$.  Since $\surfD$ is maximal, we have
\[ \Gamma \cup \Lambda \subseteq Z \cup \sigma^{-1}(Z) \cup \cdots \cup \sigma^{-(s-1)}(Z).\]

$(1)$ and $(2)$ are symmetric; we will prove $(1)$.  Let $\sh{J}$ be a nonzero graded right ideal of $\sh{S}$.
 Let $n_0 \geq s$ be such that $\sh{J}_{n_0} \neq 0$.  Let $Y$ be the subscheme of $X$ defined by $\sh{J}_{n_0}$.     
By critical transversality, there is some $n_1 \geq n_0 + s$ such that for $n \geq n_1$, we have $\sigma^{-n}(Z) \cap Y = \emptyset$.  As $\Gamma \subseteq Z \cup \sigma^{-1}(Z) \cup \cdots \cup \sigma^{-(s-1)}(Z)$, we have $\sigma^{-n}(\Gamma) \cap Y =\emptyset$ for $n \geq n_1$, as well.

 For $n \geq n_1$, let $\sh{I}_n$ be the maximal ideal sheaf on $X$ so that 
 $\sh{I}_n \supseteq \sh{J}_n$ and so that $\sh{I}_n /\sh{J}_n$ is supported on
 \[ \sigma^{-(n_1+1)} (Z) \cup \cdots \cup \sigma^{-(n+s-1)}(Z).\] 
Since 
\[\sJ_n \sS_1^{\sigma^n} = \sJ_n \sA^{\sigma^n}\sC^{\sigma^{n+1}} \subseteq \sJ_{n+1},\]
if $n \geq n_1+1$ then $\sI_n \subseteq \sI_{n+1}$.

Let $\sh{I}$ be the maximal element in the chain of the $\sh{I}_n$.  Let $m \geq n_1+1$ be such that  $\sh{I}_n = \sh{I}$ for all $n \geq m$.     Let 
\[\sh{J}' := \sh{I} \sh{D}^{\sigma^{n_1+1}} \cdots \sh{D}^{\sigma^{m-1}}.\]
Since the points in $Z$ have distinct orbits and $m > n_1$, $\sJ'$ and $\sD^{\sigma^n}$ are comaximal for $n\geq m$.
Let $n \geq m$. We must show that
\beq \label{FOO1}
\sJ_n = \sJ' \sD^{\sigma^m} \cdots \sD^{\sigma^{n-1}} \sC^{\sigma^n} = \sJ' \cap \sS_n.
\eeq

  Note first that
\beq\label{spargel} \sJ_{n_0} \sA^{\sigma^{n_0}} \sD^{\sigma^{n_0+s}} \cdots \sD^{\sigma^{n-1}} \sC^{\sigma^n} \subseteq \sJ_n 
 \subseteq \sA \sD^{\sigma^s} \cdots \sD^{\sigma^{n-1}} \sC^{\sigma^n} = \sS_n.
\eeq
If $p \in \bigcup_{j \geq m} \sigma^{-j}(Z)$, then  $p \not \in Y$.  Since the points in $Z$ have distinct orbits,  we thus have
\[ \sJ'_p = (\sJ_{n_0})_p =  \sO_{X, p}\]
and from \eqref{spargel} we obtain
\[ (\sJ_n)_p = (\sS_n)_p = (\sJ' \cap \sS_n)_p= (\sJ' \sD^{\sigma^m} \cdots \sD^{\sigma^{n-1}} \sC^{\sigma^n})_p 
 .
\]
If $p \in \bigcup_{j = n_1+1}^{m-1} \sigma^{-j}(Z)$, then (as the points in $Z$ have distinct orbits)
\[ \sI_p = (\sI_n)_p = \sO_{X,p} = (\sD^{\sigma^m} \cdots \sD^{\sigma^{n-1}} \sC^{\sigma^n})_p\]
and
\[ \sJ'_p = (\sS_n)_p = (\sJ_{n_0} \sS_{n-n_0}^{\sigma^{n_0}})_p.\]
By \eqref{spargel}, $(\sJ_n)_p = (\sS_n)_p$ and \eqref{FOO1} holds at $p$.  
If $p \in Y \cup Z \cup \cdots \cup \sigma^{-n_1}(Z)$ then by choice of $n$ we have
\[ (\sJ_n)_p = \sI_p = \sJ'_p\]
and $(\sD^{\sigma^m} \cdots \sD^{\sigma^{n-1}} \sC^{\sigma^n})_p = \struct_{X,p}$.  Thus \eqref{FOO1} holds again at $p$.  Finally, \eqref{FOO1} trivially holds at $p$ for $p \not \in Y \cup \bigcup_{j \geq 0} \sigma^{-j}(Z)$.

$(3)$. By $(1)$ and $(2)$ there are an ideal sheaf $\sh{J}$, comaximal with $\sh{D}^{\sigma^n}$ for $n \geq m$, and $\sh{K}$, comaximal with $\sh{D}^{\sigma^n}$ for $n \leq -m'$, so that
\[
 \sh{H}_n  = \sh{J} \sh{D}^{\sigma^m} \cdots \sh{D}^{\sigma^{n-1}} \sh{C}^{\sigma^n} = \sh{A} \sh{D}^{\sigma^s} \cdots \sh{D}^{\sigma^{n-m'}} \sh{K}^{\sigma^n} 
\]
for $n \gg 0$.    Let $\sh{H}' = \sh{K}^{\sigma^{m+m'}}+ \sh{J}$; note    $\sh{H}'$ is comaximal with all $\sh{D}^{\sigma^j}$.  Away from the orbits of points in $Z$, we have
$\sh{K}^{\sigma^n} = \sh{J}$ for all $n \gg 0$, and therefore for all $n$.  Therefore $\sh{H}'$ is $\sigma$-invariant. By construction $\sh{H}_n = \sh{H}' \cap \sh{S}_n = \sh{H}' \sh{S}_n$.
\end{proof}

\begin{corollary}\label{cor-ADC-noetherian}
Suppose that the tuple $\ADCD = \ADCdata$ is  transverse maximal ADC data and that  $\Lsh$ is $\sigma$-ample.  Then the ADC ring $T(\ADCD)$ and the ADC bimodule algebra $\sT (\ADCD)$ are  left and right noetherian.
\end{corollary}

\begin{proof}
Let $\sh{S} := \sT(\ADCD)$,
so that $\sh{S}_n = \sh{A} \sh{D}^{\sigma^s} \cdots \sh{D}^{\sigma^{n-1}} \sh{C}^{\sigma^n}\Lsh_n$ for $n \geq 1$.  Let $ S := T(\ADCD)$.   
 Let $\dot{\surfD} = \geomdata$ be the geometric data associated to $\surfD$; by assumption we have $\Lambda \cup \Gamma \subseteq Z\cup \cdots \cup \sigma^{-(s-1)}(Z)$.   
Since by Lemma~\ref{lem-ample1} the sequence $\{(\sh{S}_n)_{\sigma^n}\}$ is left and right ample, by Theorem~\ref{thm-VdBSerre}, to show that $S$ is noetherian it suffices to show that the bimodule algebra
$ \sh{S} $
is left and right noetherian.  By Lemma~\ref{lem-catequiv} this property does not depend on $\Lsh$, so without loss of generality we may assume that $\Lsh = \struct_X$.

By symmetry, it suffices to prove that $\sh{S}$ is right noetherian.  Let $\sh{J}$ be a graded right ideal of $\sh{S}$.  By Proposition~\ref{prop-circ2}(1), there are an ideal sheaf $\sh{J}'$ on $X$ and an integer $m \geq s$ such that for $n \geq m$, 
\[ \sh{J}_n = \sh{J}' \sh{D}^{\sigma^{m}} \cdots \sh{D}^{\sigma^{n-1}} \sh{C}^{\sigma^n}\]
and $\sJ'$ is comaximal with $\sD^{\sigma^n}$.
We claim that $\sh{J}$ is generated by $\sh{J}_{\leq m+2s}$. 

 This is a straightforward computation.  
 Since $\mc{AC} \subseteq \sD\cdots \sD^{\sigma^{s-1}}  \subseteq \sA \cap \sC$ and the points of $Z$ have distinct orbits, note that 
 \[\sh{AC}\sh{D}^{\sigma^s}\cdots \sh{D}^{\sigma^{2s-1}} + \sh{D}\sh{D}^{\sigma} \cdots \sh{D}^{\sigma^{s-1}}(\sh{AC})^{\sigma^s} = \sh{D}\cdots \sh{D}^{\sigma^{2s-1}}.\]
  Let $k \geq 2s$. We have
\begin{multline*}
\sh{J}_m \sh{S}_k^{\sigma^m} + \sh{J}_{m+s}\sh{S}_{k-s}^{\sigma^{m+s}} = \\
\sh{J}' \sh{C}^{\sigma^m} \sh{A}^{\sigma^m} \sh{D}^{\sigma^{m+s}} \cdots \sh{D}^{\sigma^{m+k-1}}\sh{C}^{\sigma^{m+k}} \\ +
\sh{J}' \sh{D}^{\sigma^m} \cdots \sh{D}^{\sigma^{m+s-1}} \sh{C}^{\sigma^{m+s}} \sh{A}^{\sigma^{m+s}} 
\sh{D}^{\sigma^{m+2s}} \cdots \sh{D}^{\sigma^{m+k-1}} \sh{C}^{\sigma^{m+k}}  \\
= \sh{J}'  \sh{D}^{\sigma^{m}}  \cdots  \sh{D}^{\sigma^{m+k-1}} \sh{C}^{\sigma^{m+k}}  = \sh{J}_{m+k}.
 \end{multline*}
  Thus $\sh{J}_{\geq m+2s} = (\sh{J}_{m} \sh{S} + \sh{J}_{m+s} \sh{S})_{\geq m+2s}$.  The claim follows, and $\sh{J}$ is coherent.

\end{proof}

Recall the notation that if $I$ is a right ideal of a ring $R$, then 
\[\I^{\ell}_R(I) = \{ r \in R \st rI \subseteq I\}\]
 is the maximal subring of $R$ so that $I$ is a two-sided ideal.  
We call this ring the {\em left  idealizer of $I$ in $R$}.  If $I$ is a left ideal of $R$, we similarly define the {\em right 
 idealizer} of $I$ in $R$ to be
\[ \I^{r}_R(I) = \{ r \in R \st Ir \subseteq I\}.\] 
We similarly define $\I^r_{\sh{R}}(\sh{I})$, respectively $\I^{\ell}_{\sR}(\sI)$, for a left, respectively right,  ideal $\sh{I}$ of a bimodule algebra $\sh{R}$.

We  recall the criteria for an idealizer bimodule algebra to be left or right noetherian.

\begin{proposition}\label{prop-ibma-lnoeth}
{ \em (\cite[Proposition~3.3.3]{S-thesis}, cf. \cite[Proposition~2.2]{R-idealizer})}
Let $X$ be a projective variety, and let $\sigma \in \Aut(X)$. Let $\sh{B}$ be a noetherian graded $(\struct_X, \sigma)$-sub-bimodule algebra of $\sh{B}(X, \struct_X, \sigma)$, and let $\sh{I} = \bigoplus (\sh{I}_n)_{\sigma^n}$ be  a graded right ideal of $\sh{B}$. Let $\sh{R} := \I^{\ell}_{\sh{B}}(\sh{I})$. 
 Suppose that $\sh{R}_n = \sh{I}_n$  for all $n \gg 0$.    Then $\sh{R}$ is left noetherian if and only if for all graded left ideals $\sh{J}$ of $\sh{B}$ we have
\[ (\sh{I}  \cap \sh{J})_n =( \sh{IJ})_n\]
for $n \gg 0$. \qed
\end{proposition}

\begin{proposition}\label{prop-ibma-rnoeth}
{\em (\cite[Lemma~3.9]{S-idealizer}, cf. \cite[Lemma~1.2]{Stafford1985}, \cite[Proposition~2.1]{R-idealizer})}
Let $X$ be a projective variety, and let $\sigma \in \Aut(X)$.  Let $\sh{B}$ be  a right noetherian graded $(\struct_X, \sigma)$-sub-bimodule algebra of the twisted bimodule algebra 
$\sh{B}(X, \struct_X, \sigma),$
 and let $\sh{I} = \bigoplus (\sh{I}_n)_{\sigma^n} $ be  a nonzero graded right ideal of $\sh{B}$.  Let $\sh{R} := \I^{\ell}_{\sh{B}}(\sh{I})$. 
 Suppose that for all graded right ideals $ \sh{J} \supseteq \sh{I}$ 
of $\sh{B}$,  for $n \gg 0$ we have
\[ \sh{B}_n \cap \bigcap_{m \geq 0} (\sh{J}_{n+m}: \sh{I}_m^{\sigma^n}) = \sh{J}_n.\]
Then $\sh{R}$ is right noetherian.
\qed
\end{proposition}

It is straightforward to obtain defining data for left and right ideals of idealizer bimodule algebras, and we do this next.

\begin{lemma}\label{lem-circ1}
Let $X$ be a projective variety, let $\sigma \in \Aut(X)$, and let $\Lsh$ be an invertible sheaf on $X$.  Let 
\[\sh{S} := \bigoplus_{n \geq 0} (\sh{S}_n)_{\sigma^n}\]
 be a noetherian sub-bimodule algebra of $\sh{B}(X, \Lsh, \sigma)$, and let $\sh{I} = \bigoplus (\sh{I}_n)_{\sigma^n}$ be a graded right ideal of $\sh{S}$.  Let $\sh{R} := \I^{\ell}_{\sh{S}}(\sh{I})$, 
and assume that $\sh{R}$ is also noetherian and that $\sh{R}_n = \sh{I}_n$ for $n \gg 0$.  

$(1)$
Let $\sh{J} = \bigoplus (\sh{J}_n)_{\sigma^n}$ be a graded right ideal of $\sh{R}$.  Then there is a right ideal $\sh{J}' \subseteq \sh{I}$ of $\sh{S}$ such that 
\[ \sh{J}_n = (\sh{J}')_n\]
for $n \gg 0$.

$(2)$ Let $\sh{K}$ be a graded left ideal of $\sh{R}$.  Then there is a graded left ideal $\sh{K}'$ of $\sh{S}$ 
such that
\[ \sh{K}_n = (\sh{I} \cap \sh{K'} )_n = (\sh{IK}')_n\]
for $n \gg 0$.

$(3)$  Let $\sh{H}$ be a graded ideal of $\sh{R}$.  Then there is a graded ideal $\sh{H}'$ of $\sh{S}$ so that
$\sh{H}$, $\sh{I} \cap \sh{H}'$, and $\sh{IH}'$ are equal in large degree.  
\end{lemma}
\begin{proof}
$(1)$.  Fix $\sh{J}$.  Since $\sh{R}$ is noetherian, there is an integer $k$ such that $\sh{J}$ is generated in degree $\leq k$.  Let $\sh{J'} := \sh{JI}$.  Then $\sh{J}'$ is a right ideal of $\sh{S}$.   Since $\sh{R}_n = \sh{I}_n$ for $n \gg 0$, we have
\[\sh{J}'_n = (\sh{J I})_n = (\sh{J}_{\leq k} \sh{R})_n = \sh{J}_n\]
for $n \gg k$.

$(2)$. Fix $\sh{K}$ and let $\sh{K'} := \sh{SK}$.  A similar argument shows that for $n \gg 0$ that $(\sh{I}  \sh{K}')_n = \sh{K}_n$. 
For $n \gg 0$ we have  $(\sh{I} \cap \sh{K'})_n = (\sh{IK}')_n$, by  Proposition~\ref{prop-ibma-lnoeth}.

$(3)$.  The construction in part $(1)$ shows that if we replace $\sh{H}$ by $\sh{H}_{\geq n}$ for some $n\gg 0$, we may assume without loss of generality that $\sh{H}= \sh{HS}$ is  a right ideal of $\sh{S}$.  The proof of $(2)$ shows that
\[ \sh{H} = \sh{I} \cap \sh{SH} = \sh{ISH}\]
in large degree.  Let $\sh{H}' := \sh{SH}$. 
\end{proof}

The computations in the next lemma will allow us to apply Proposition~\ref{prop-ibma-lnoeth} to show  for arbitrary tranverse ADC data $\surfD$ that  the bimodule algebra $\sT(\surfD)$ is noetherian.  

\begin{lemma}\label{lem-easy2}
Let 
\[\surfD = \surfdata \]
be  transverse ADC data, and let 
\[ \sh{T} := \sh{T}(\surfD).\]
Let $\sh{J}$ be an ideal sheaf on $X$. Let $\sh{A}' \supset  \sh{A}$ and $\sh{C}' \supset \sh{C}$ be any pair maximal with respect to 
\[ \sh{A}'\sh{C}' \subseteq \sD \sD^{\sigma} \cdots \sD^{\sigma^{s-1}}.\]

$(1)$  Let 
\[ \sh{R} := \sh{T}(X, \Lsh, \sigma, \sh{A}, \sh{D}, \sh{C}',s).\]
Let $\sh{K}$ be the left ideal of $\sh{R}$ defined by
\[ \sh{K}_n := \sh{R}_n \cap \sh{J}^{\sigma^n}\Lsh_n.\]
 Then for $n \gg 0$, we have
\[
\sh{K}_n = \sh{R}_n \cap \bigcap_{m\geq 1} (\sh{K}_{n+m}: \sh{T}_m)^{\sigma^{-m}}.\]

$(2)$  Let
\[ \sh{R}' := \sh{T}(X, \Lsh, \sigma, \sh{A}', \sh{D}, \sh{C},s).\]
Let $\sh{K}'$ be the right ideal of $\sh{R}'$ defined by
\[ \sh{K}'_n := \sh{R}'_n \cap \sh{J}\Lsh_n.\]
Then for $n \gg 0$, we have
\[ \sh{K}'_n = \sh{R}'_n \cap \bigcap_{m \geq 1} (\sh{K}'_{n+m}: \sh{T}_m^{\sigma^n}).\]
\end{lemma}
\begin{proof}
Without loss of generality we may assume that $\Lsh = \struct_X$.  
$(1)$ and $(2)$ are similar; we prove $(1)$.  Since
\[ \sh{T}_m \sh{K}_n^{\sigma^m} \subseteq \sh{R}_{n+m} \cap \sh{J}^{\sigma^{n+m}} 
= \sK_{n+m}\]
the inclusion $\subseteq$ is trivial.  

  For the other inclusion, we certainly have
\begin{multline*}
\bigcap_{m \geq 1}(\sh{K}_{n+m}: \sh{T}_m)^{\sigma^{-m}} \cap \sh{R}_n = 
\bigcap_{m \geq 1}\bigl( (\sh{R}_{n+m} \cap \sh{J}^{\sigma^{n+m}}): \sh{T}_m\bigr)^{\sigma^{-m}} \cap \sh{R}_n \\
\subseteq
\bigcap_{m \geq 1}(\sh{J}^{\sigma^{n+m}}: \sh{T}_m)^{\sigma^{-m}} \cap \sh{R}_n 
= \bigcap_{m \geq 1} (\sJ^{\sigma^n}: \sT_m^{\sigma^{-m}}) \cap \sR_n.
\end{multline*}

Let
\[ \sJ = \sK_1 \cap \cdots \cap \sK_{\ell}\]
be a minimal primary decomposition of $\sJ$, where $\sK_i$ is $\sQ_i$-primary. Since $\sD$ and $\sC$ are cosupported at points of infinite order, we may  choose $n_0$ so that for $n \geq n_0$, neither $\sC$ or any $\sD^{\sigma^j}$, where $j<0$, are contained in any $\sQ_i^{\sigma^n}$.   By \cite[Lemma~2.13(1)]{S-idealizer}, for $n \geq n_0$
\[ (\sJ^{\sigma^n}: \sT_m^{\sigma^{-m}}) \subseteq (\sJ^{\sigma^n}: \sA^{\sigma^{-m}}).\]

Fix $n \geq n_0$.  Transversality of $\surfD$ implies that no primary component of $\sA^{\sigma^{-m}}$ is contained in any $\sQ_i^{\sigma^n}$ for $m \gg 0$.  Thus by  \cite[Lemma~2.13(2)]{S-idealizer},
\[(\sJ^{\sigma^n}: \sA^{\sigma^{-m}}) = \sJ^{\sigma^n}\]
 for $m \gg 0$.  Thus
\[ 
 \sR_n \cap \bigcap_{m \geq 1}(\sJ^{\sigma^n}: \sA^{\sigma^{-m}}) = \sR_n \cap \sJ^{\sigma^n} = \sK_n,
\]
and (1) holds.
\end{proof}

\begin{corollary}\label{cor-easy3}
Let
\[\surfD = \surfdata\]
 be  transverse ADC data. 
Let 
\[ \sh{C}' := (\sD \sD^{\sigma} \cdots \sD^{\sigma^{s-1}}:  \sh{A})\]
and let
\[ \sh{A}':= (\sD \sD^{\sigma} \cdots \sD^{\sigma^{s-1}} : \sh{C}).\]
Let 
\[ \mb{E} := (X, \Lsh, \sigma, \sh{A}, \sh{D}, \sh{C}',s)\]
and let
\[ \mb{E}' := (X, \Lsh, \sigma, \sh{A}', \sh{D}, \sh{C},s).\]
Let
$\sh{T} := \sh{T}(\surfD)$, let $\sh{R}:= \sh{T}(\mb{E})$, and let
$\sh{R}' := \sh{T}(\mb{E}')$.  
 Then
$ \sh{T}$, $ \I^{r}_{\sh{R}}(\sh{T}_{\geq 1})$, and $\I^{\ell}_{\sh{R}'}(\sh{T}_{\geq 1})$ are all equal in large degree. 
\end{corollary}
\begin{proof}
Without loss of generality, we may assume that  $\Lsh = \struct_X$.  

Note that  $\sT_{\geq 1}$ is a right ideal of $\sR'$ and a left ideal of $\sR$.  Since $\surfD$ is transverse, for some $k\geq s$ we have
\[ \sT_n = (\sA \sD^{\sigma^s} \cdots \sD^{\sigma^{k-1}}) \cap \sR'_n\]
for $n \geq k$.   Put
\[ \sJ := \sA \sD^{\sigma^s} \cdots \sD^{\sigma^{k-1}}\]
and let $\sK' = \bigoplus_n \left( \sR'_n \cap \sJ \right)$.  Then $\sK'$ and $\sT$ are equal in large degree, and by Lemma~\ref{lem-easy2}(2)
\[ 
 \sT_n = \sK'_n = \sR'_n \cap \bigcap_{m \geq 1}(\sK'_{n+m}:\sT_m^{\sigma^n}) 
= \sR'_n \cap \bigcap_{m\geq 1}(\sT_{n+m}:\sT_{m}^{\sigma^n}) = \Bigl( \I^{\ell}_{\sR'}(\sT_{\geq 1}) \Bigr)_n  
\]
for $n \gg 0$.

The proof that $\sT$ and $\I_{\sR}^{r}(\sT_{\geq 1})$  are equal in large degree is symmetric.
\end{proof}

We are now ready to show that if $\surfD$ is transverse ADC data, then $T(\surfD)$ is noetherian.

\begin{proposition}\label{prop-ringsnoeth}
Suppose that  $\surfD = \surfdata$ is transverse   ADC data  and that $\Lsh$ is $\sigma$-ample.  
Let 
$\sh{T} := \sh{T} (\surfD)$
 and let 
$T := T (\surfD)$.
   Then both $\sh{T}$ and $T$ are noetherian.
\end{proposition}
\begin{proof}
  By Lemma~\ref{lem-ample1}(2), the sequence of bimodules $\{ (\sh{T}_n)_{\sigma^n} \}$ is left and right ample.  Thus by Theorem~\ref{thm-VdBSerre}, it suffices to prove that $\sh{T}$ is right and left noetherian.  Without loss of generality we may assume that $\Lsh = \struct_X$.  

If  $\surfD$ is maximal, this is Corollary~\ref{cor-ADC-noetherian}.   Suppose that 
$\surfD$ is right maximal but not left maximal.  Let 
\[ \sh{A}' := (\sD \sD^{\sigma} \cdots \sD^{\sigma^{s-1}} : \sh{C})\]
and let
\[ \sh{S} := \sT(X, \struct_X, \sigma, \sh{A}', \sh{D}  , \sh{C}).\]
By Corollary~\ref{cor-ADC-noetherian}, $\sh{S}$ is left and right noetherian.

Now, $\sh{T}_{\geq 1}$ is a graded right ideal of $\sh{S}$.  Let 
$\sh{K} \supseteq \sh{T}_{\geq 1}$ be another graded right ideal of $\sh{S}$.  By Proposition~\ref{prop-circ2},  there are an ideal sheaf $\sh{J}$ on $X$ and an integer $k \geq 0$ so that for $n \geq k$ we have
\[ \sh{K}_n = \sh{J} \cap \sh{S}_n.\]
Let
\[ \sh{F}_n := \sh{S}_n \cap \bigcap_{m \geq 0}(\sh{K}_{n+m}: \sh{T}_m^{\sigma^n}).\]
By Lemma~\ref{lem-easy2},  $\sh{F}_n = \sh{K}_n$ for $n \gg 0$. 
By Corollary~\ref{cor-ADC-noetherian}, $\sh{S}$ is noetherian. Thus by Proposition~\ref{prop-ibma-rnoeth},
$\I^{\ell}_{\sh{S}}(\sh{T}_{\geq 1})$ is right noetherian.  By Corollary~\ref{cor-easy3}, $\sh{T}$ is also right noetherian.

Now suppose that $\sh{K}$ is a graded left ideal of $\sh{S}$; by Proposition~\ref{prop-circ2}, there are  an ideal sheaf $\sh{J}$ and  on $X$  and an integer $k \geq s $ so that for $n \geq k+s$ we have:
\[\sh{K}_n = \sh{A}' \sh{D}^{\sigma^s} \cdots \sh{D}^{\sigma^{n-k}} \sh{J}^{\sigma^n},\]
and 
$\sh{J}$ and $\sh{D}^{\sigma^j}$ are comaximal for $j \leq -k$.
Then for $n >  k +s$, we have 
\beq\label{hillary}
(\sh{T} \cap \sh{K})_n  = \sh{A} \sh{D}^{\sigma^s} \cdots \sh{D}^{\sigma^{n-1}} \sh{C}^{\sigma^n} \cap
\sh{A}' \sh{D}^{\sigma^s} \cdots \sh{D}^{\sigma^{n-k}} \sh{J}^{\sigma^n}.
\eeq
Transversality of the defining data for $\sh{T}$ implies that 
\[ \sh{A} \cap \sh{J}^{\sigma^n}  = \sh{A}  \sh{J}^{\sigma^n}\]
for $n \gg 0$.  
Thus \eqref{hillary} is equal to 
\[ \sh{A} \sh{D}^{\sigma^s} \cdots \sh{D}^{\sigma^{n-k} }\sh{J}^{\sigma^n}\]
for $n \gg 0$.

On  the other hand, for $n \geq 2k+2s$ we have
\[ ((\sh{T}_{\geq 1})\cdot \sh{K})_n \supseteq \sh{T}_k (\sh{K}_{n-k})^{\sigma^k} + \sh{T}_{k+s} (\sh{K}_{n-k-s})^{\sigma^{k+s}}
  = \sh{A} \sh{D}^{\sigma^s} \cdots \sh{D}^{\sigma^{n-k}} \sh{J}^{\sigma^n}.\]
  Thus 
  \[ ((\sh{T}_{\geq 1})\sh{K})_n \supseteq (\sh{T} \cap \sh{K})_n\]
  for $n \gg 0$.  As the other containment is automatic, by Proposition~\ref{prop-ibma-lnoeth} both $\I^{\ell}_{\sS}(\sT_{\geq 1})$ and $\sh{T}$ are left noetherian.

We now consider the general case.  Given  transverse ADC data 
\[\surfD = (X, \struct_X, \sigma, \sh{A}, \sh{D}, \sh{C}),\]
 let 
$\sh{C}':= (\sD \sD^{\sigma} \cdots \sD^{\sigma^{s-1}}: \sh{A})$ and let
\[ \mathbb{F} := (X, \struct_X, \sigma, \sh{A}, \sh{D}, \sh{C}').\]
Let 
\[\sh{R} := 
\sh{T}(\mathbb{F}).\]
Sincd $\mb{F}$ is right maximal,  $\sh{R}$ is noetherian.  

Note $\sT_{\geq 1}$ is a left ideal of $\sR$.  
Let $\sh{K} \supseteq \sh{T}_{\geq 1}$ be a graded 
left ideal of $\sh{R}$.  By Proposition~\ref{prop-circ2} and Lemma~\ref{lem-circ1}(2), there is an ideal sheaf $\sh{J}$ on $X$ so that 
\[ \sh{K}_n = \sh{J}^{\sigma^n} \cap \sh{R}_n\]
for $n \gg 0$.  By Lemma~\ref{lem-easy2}(1), we have
\[ \sh{K}_n = \sh{R}_n \cap \bigcap_{m \geq 1} (\sh{K}_{n+m} :\sh{T}_m)^{\sigma^{-m}}  \]
for $n \gg 0$.  By Corollary~\ref{cor-easy3}, $\sh{T}$ and $\I_{\sh{R}}^{r}(\sh{T}_{\geq 1})$ 
are equal in large degree.  By the left-handed version of Proposition~\ref{prop-ibma-rnoeth}, $\sh{T}$ is left noetherian.  

By symmetry, $\sh{T}$ is also right noetherian.  
\end{proof}

We record a result on two-sided ideals of the rings $T(\surfD)$, which will be useful later in this paper and in \cite{S-surfclass}.

\begin{proposition}\label{prop-ideals}
Let $\surfD= \surfdata$ be transverse ADC data, where $\Lsh$ is $\sigma$-ample.  Let $\sh{T} := \sh{T}(\surfD)$ and let $T:= T(\surfD)$.  Let $K$ be a two-sided graded ideal of $T$.  Then there is a $\sigma$-invariant ideal sheaf $\sh{K}$ on $X$ so that 
\[ K_n = H^0(X, \sh{KT}_n) = H^0(X, \sh{KL}_n \cap \sh{T}_n)\]
for $n \gg 0$.
\end{proposition}

\begin{proof}
Let 
\[ \sh{A}':=(\sD \sD^{\sigma} \cdots \sD^{\sigma^{s-1}}: \sh{C})\]
and let
\[ \sh{C}':= (\sD \sD^{\sigma} \cdots \sD^{\sigma^{s-1}}: \sh{A}').\]
Let $\mb{E}:= (X, \Lsh, \sigma, \sh{A}', \sh{D}, \sh{C},s)$ and let $\mb{F}:= (X, \Lsh, \sigma, \sh{A}', \sh{D}, \sh{C}',s)$.  Let $\sh{R}:= \sh{T}(\mb{E})$ and let $\sh{S}:=\sh{T}(\mb{F})$.
Note that both $\mb{F}$ and $\mb{E}$ are transverse, and that $\mb{F}$ is maximal.  By Proposition~\ref{prop-ringsnoeth}, $\sh{T}$,  $\sh{R}$, and $\sh{S}$ are noetherian.  By Corollary~\ref{cor-easy3}, $\sh{T}$ is equal in large degree to a left 
 idealizer inside $\sh{R}$, and $\sh{R}$ is equal in large degree to a right 
 idealizer inside $\sh{S}$.

By Theorem~\ref{thm-VdBSerre} there is a two-sided ideal $\sh{F} = \bigoplus (\sh{F}_n)_{\sigma^n}$ of $\sh{T}$ so that $K_n = H^0(X, \sh{F}_n)$ for $n \gg 0$.
Applying Lemma~\ref{lem-circ1}(3) twice, we obtain a two-sided ideal $\sh{H}$ of $\sh{S}$ so that $\sh{F}_n = \sh{H}_n \cap \sh{T}_n$ for $n \gg 0$.  By Proposition~\ref{prop-circ2}(3), there is a $\sigma$-invariant ideal sheaf $\sh{K}$ on $X$ so that 
\[ \sh{H}_n = \sh{KS}_n = \sh{KL}_n \cap \sh{S}_n\]
for $n \gg 0$.  Thus $\sh{F}_n = \sh{KL}_n \cap \sh{T}_n$ for $n \gg 0$.  Transversality of $\surfD$ and $\sigma$-invariance of $\sh{K}$ imply that
\[ \sh{KL}_n \cap \sh{T}_n = \sh{KT}_n.\]
\end{proof}

We now prove the converse to Proposition~\ref{prop-ringsnoeth}.  We do this in several steps.

\begin{lemma}\label{lem-notrighttransverse}
Let $\surfD = \surfdata$ be ADC data, where $\Lsh$ is $\sigma$-ample.  Let $T:=T(\surfD)$ and let $\sh{T}:=\sh{T}(\surfD)$.  Let $Z$ be the cosupport of $\sh{D}$ and let $\Gamma$ be the cosupport of $\sh{C}$.  If $T$ or $\sh{T}$ is right noetherian, then the sets $\{ \sigma^n Z\}_{n \leq 0}$ and $\{\sigma^n\Gamma\}_{n \leq 0}$ are critically transverse. 
\end{lemma}
\begin{proof}
 Suppose that $T$ is right noetherian and that critical transversality fails.  Let 
\[ \sh{A}':= (\sD \sD^{\sigma} \cdots \sD^{\sigma^{s-1}}: \sh{C})\]
and let 
\[S := T(X, \Lsh, \sigma, \sh{A}', \sh{D}, \sh{C},s).\]
Then $T_{\geq 1}$ is a right ideal of $S$.  It is easy to see that $S_T$ is isomorphic to a right ideal of $T$ and is thus finitely generated.  Thus $S$ is right noetherian and we may assume without loss of generality that $\surfD$ is left maximal.

We claim that  there are a subscheme $Y$ of $X$, a point $p \in Z \cup \Gamma$, and an integer $n_0$ so that
\begin{enumerate}
\item
$\sh{I}_Y\Lsh_n \cap \sh{T}_n$ is globally generated for $n \geq n_0$; and 
\item  the set $\{ \sigma^n (p)\}_{n \leq 0} \cap Y$ is infinite. 
\end{enumerate} 
Assume this for the moment, and suppose that $p \in Z$.  Let $\sh{F}_n:= \sh{I}_Y \Lsh_n \cap \sh{T}_n$ and let 
\[ F:= \bigoplus_{n \geq 0} H^0(X, \sh{F}_n).\]
Note $F$ is a right ideal of $T$.  
Fix $k \in \NN$ and let $n \geq k +s$.  Let $m \geq n, n_0$ be so that $\sigma^{-(m-1)}(p) \in Y$.  

Let $\sh{E}$ be the ideal sheaf so that $\sh{E}_p = (\sh{DC}^{\sigma})_p$ and so that $\sh{E}$ is cosupported at $p$.  
Since $\sigma^{-(m-1)}(p) \in Y$, we have
\[ (\sh{F}_m )_{\sigma^{-(m-1)}(p)} = \bigl( (\sh{I}_Y \cap \sh{E}^{\sigma^{m-1}} )\Lsh_m\bigr)_{\sigma^{-(m-1)}(p)} \supsetneqq (\sh{I}_Y \sh{E}^{\sigma^{m-1}} \Lsh_m)_{\sigma^{-(m-1)}(p)}.\]
Now,  
\[(F_{\leq k} \cdot T)_m \subseteq
H^0(X, \sh{F}_1 \sh{T}_{m-1}^{\sigma} + \cdots + \sh{F}_k \sh{T}_{m-k}^{\sigma^k})
\subseteq H^0(X, \sh{I}_Y \sh{E}^{\sigma^{m-1}} \Lsh_m \cap \sh{T}_m).\]
Since $\sh{F}_m$ is globally generated, $(F_{\leq k} \cdot T)_m$  is strictly contained in $F_m$.  Thus $F$ is not finitely generated as a right ideal.

If the claim holds and  $p \in \Gamma$, the argument that $F$ is not finitely generated is similar.  

Thus it suffices to prove the claim.  
Recall that $\sA \supseteq \sD \sD^{\sigma} \cdots \sD^{\sigma^{s-1}}$ because $\surfD$ is left maximal. 
If the orbits of all points in $Z \cup \Gamma$ are Zariski-dense, then the sequence $\{(\sh{T}_n)_{\sigma^n}\}$ is left and right ample.  Let $Y$ be any subscheme that has infinite intersection with some $\{\sigma^n(p)\}_{n \leq 0}$ for $p \in \Gamma \cup Z$.  By ampleness, $\sh{I}_Y \Lsh_n \cap \sh{T}_n$ is globally generated for $n \gg 0$.

Thus suppose that there is some point in $Z \cup \Gamma$ whose orbit is not dense.  
Let $\sh{I}_n := \sh{T}_n \Lsh_n^{-1}$.  Write 
\[ \sh{I}_n = \sh{J}_n \cap \sh{K}_n,\]
where $\sh{J}_n$ is cosupported at points whose orbits are not dense, and $\sh{K}_n$ is cosupported at points with dense orbits.    Let $W_n$ be the subscheme defined by $\sh{J}_n$.  Let $Y$ be the Zariski closure of the schemes $\{\sigma^m W_n\}_{n \in \NN, m \in \ZZ}$.   Note $Y$ is a proper $\sigma$-invariant subscheme of $X$.  
  By Lemma~\ref{lem-ample1}, the sequence of bimodules $\{ (\sh{K}_n \Lsh_n)_{\sigma^n} \}$ is left and right ample, so $\sh{I}_Y \sh{K}_n \Lsh_n$ is globally generated for $n \gg 0$.  But $\sh{I}_Y \subseteq \sh{J}_n$ for all $n$, and $\sh{I}_Y $ is comaximal with $\sh{K}_n$.  Thus $\sh{I}_Y \sh{K}_n \Lsh_n = \sh{I}_Y \Lsh_n \cap \sh{T}_n$ and the claim holds. 
  
  If $\sh{T}$ is right noetherian and  transversality fails, then as above we may assume without loss of generality that $\surfD$ is left maximal.  Let $Y$ and $\sh{F}_n$ be as above.  The proof above shows that $\sh{F}:= \bigoplus \sh{F}_n$ is a right ideal of $\sh{T}$ that  is not coherent.
\end{proof}

Recall that if $X$ is a projective variety, $\sigma \in \Aut(X)$, $\Lsh$ is a $\sigma$-ample invertible sheaf, and $Z$ is a 0-dimensional subscheme of $X$ supported on infinite $\sigma$-orbits, then the {\em na\"ive blowup} 
\[ S = S(X, \Lsh, \sigma, Z)\]
is defined as
\[ S = \bigoplus_{n \geq 0} H^0(X, \sh{I}_n \Lsh_n),\]
where
\[ \sh{I}_n := \sh{I}_Z \sh{I}_Z^{\sigma} \cdots \sh{I}_Z^{\sigma^{n-1}}.\]
We define
\[ \sh{S}(X, \Lsh, \sigma, Z):= \bigoplus_{n \geq 0} (\sh{I}_n \Lsh_n)_{\sigma^n}.\]
By \cite[Theorem~3.1]{RS-0}, if all points in $Z$ have critically dense orbits, then $S$ is noetherian.  
We note that we can  prove the converse   using similar methods as in the proof of Lemma~\ref{lem-notrighttransverse}.
\begin{proposition}\label{cor-NB}
Let $X$ be a projective variety of dimension $ \geq 2$, let  $\sigma \in \Aut (X)$, and let $\Lsh$ be a $\sigma$-ample invertible sheaf on $X$.  Let $Z$ be a 0-dimensional subscheme of $X$ supported at points of infinite order, and let $S$ be the na\"ive blowup algebra $S:= S(X, \Lsh, \sigma, Z)$.  Then $S$ is noetherian if and only if all points in $Z$ have critically dense orbits.
\end{proposition}
\begin{proof}
If all points have critically dense orbits, then $S$ is noetherian by \cite[Theorem~3.1]{RS-0}.
If all points in $Z$ have dense orbits but some orbit is not critically dense, then $S$ is not noetherian by  \cite[Proposition~3.16]{RS-0}.  Thus it suffices to suppose that some point in $Z$ has a non-dense orbit, and show that $S$ is not noetherian.  This follows as in the proof of Lemma~\ref{lem-notrighttransverse}.
\end{proof}

We now show that if $\Lsh$ is $\sigma$-ample, then transversality of ADC data characterizes when the algebras $T(\surfD)$ are noetherian.

\begin{theorem}\label{thm-noeth}
Let $\surfD = \surfdata$ be ADC data, where $\Lsh$ is $\sigma$-ample.  Let $T:=T(\surfD)$ and let $\sh{T} := \sh{T}(\surfD)$.  Then the following are equivalent:

$(1) $ $T$ is noetherian;

$(2)$ $\sh{T}$ is noetherian;

$(3)$  $\surfD$ is transverse.
\end{theorem}
\begin{proof}
$(3)$ $\implies$ $(1)$, $(2)$ is  Proposition~\ref{prop-ringsnoeth}. 

$(1)$ $\implies $ $(3)$, $(2) \implies (3)$.  Suppose that  $T$ or $\sh{T}$ is noetherian.  
Let 
\[\dot{\surfD} = \geomdata\]
 be the geometric data associated to $\surfD$.  
  By Lemma~\ref{lem-notrighttransverse}, the sets $\{\sigma^n Z\}_{n \leq 0}$ and $\{\sigma^n \Gamma\}_{n \leq 0}$ are critically transverse.  In particular, $Z$ and $\Gamma$ are supported on dense orbits, and $\{(\sh{T}_n)_{\sigma^n}\}$ is a right ample sequence, by Lemma~\ref{lem-ample1}.  

Our next step is to show that $\Omega$ is locally principal and contains no points or curves of finite order.  
Let $W$ be the subscheme defined by $\sA$.     

We claim that for any $\sigma$-invariant subscheme $Y$, we have $\shTor_1^X(\struct_{W}, \struct_Y) = 0$; that is, that $\sh{A} \cap \sh{I}_Y = \sh{A} \sh{I}_Y$.  Suppose this fails for some $\sigma$-invariant $Y$, so 
\[ \sh{A} \cap \sh{I}_Y \supsetneqq \sh{A}\sh{I}_Y.\]
Let 
\[ \sh{J}:= \bigoplus_{n \geq 0}  \sh{I}_Y \Lsh_n \cap \sh{T}_n\] 
and let $J:= H^0(X, \sh{J})$.   Right ampleness and Corollary~\ref{cor-FOOBAR} imply that  there exists a $k$ so that $J_n$ generates
$\sh{J}_n$ for $n \geq k$.  
Fix $k' \geq k$ and consider the left ideal $T (J_{\leq k'})$.   For $n \geq k'+s$ we have
\[ T(J_{\leq k'})_n \subseteq H^0(X, \sI_Y \sh{A} \sh{D}^{\sigma^s}\cdots \sh{D}^{\sigma^{n-1}} \sh{C}^{\sigma^n} \Lsh_n).\]  
This is not equal to $J_n$ by global generation of $\sh{J}$ and by assumption on $Y$.  Thus ${}_TJ$ is not finitely generated.  Likewise, ${}_{\sh{T}} \sh{J}$ is not coherent.  Since $T$ or $\sh{T}$ is noetherian,  no such $Y$ can exist.  

In particular, $W$ does not contain any points or components of finite order (and thus neither does $\Omega$). 
We show that $\Omega$ is locally principal; we only need to check this at singular points of $X$.  
Suppose then that $W \cap X^{\sing} \neq \emptyset$.  Since $W$ contains no points or components of finite order, we may assume that $X^{\sing}$ is a curve and $W \cap X^{\sing}$ is 0-dimensional; further,  $X^{\sing}$ is smooth at all  points of $W \cap X^{\sing}$.  Moreover, since $X$ fails Serre's condition $S_2$ at only finitely many points, $X$ is $S_2$ at all points of $X^{\sing} \cap W \supseteq X^{\sing} \cap \Omega$.  

 Let $x \in X^{\sing} \cap \Omega$.  Let $K \subseteq \sO_{X,x}$ be the ideal defining $X^{\sing}$ at $x$.  Then $(\sA_x+ K)/K$ defines $W \cap X^{\sing} \subset X^{\sing}$.  Since $O_{X,x}/K$ is a regular local ring of dimension 1, $\sA_x+K$ is principal modulo $K$, and there is some $f \in \sA_x$ so that $\sA_x+K = (f) +K$.  As in \cite[Lemma~7.7]{S-idealizer}, $\sA_x$ is a principal ideal of $\sO_{X,x}$.  Since $\sO_{X, x}$ satisfies $S_2$, the associated primes of $(f)$ are all height 1.  By definition of $\Omega$, this says precisely that $\sA_x = (\sI_{\Omega})_x$.  Thus $\Omega$ is locally principal at $x$.

As in the comments after Proposition~\ref{prop-transverse-char0}, for appropriate $\surfD'$ we have  $T \cong T(\surfD')^{op}$ and $\sh{T} \cong \sh{T}(\surfD')^{op}$, and by symmetry   $\{\sigma^n Z\}_{n \geq 0}$ and $\{\sigma^n \Lambda\}_{n \geq 0}$ are also critically transverse.  

By Lemma~\ref{lem-CT}, it suffices to prove that $\Omega$ meets orbits of points finitely often, and by symmetry it suffices to show that $\Omega$ meets backward orbits finitely often.  Let $p$ be a point (necessarily of infinite order) so that $\{\sigma^n (p)\}_{n \leq 0} \cap \Omega$ is infinite.  We may assume that 
\[p \not\in \Gamma \cup \{\sigma^nZ\}_{n \geq 0} \cup \{\sigma^n \Lambda\}_{n \geq 0}.\]
Let 
\[\sh{A}':=  (\sD \sD^{\sigma} \cdots \sD^{\sigma^{s-1}}:\sh{C}) \cap \sI_{\Lambda}.\]
Note 
\[\Cosupp \sh{A}' \subseteq Z \cup \sigma^{-1}(Z) \cup \cdots \cup \sigma^{-(s-1)}(Z) \cup \Lambda.\]
  Let
$\ADCD := (X, \Lsh, \sigma, \sh{A}', \sh{D}, \sh{C},s)$.  Let $S:=T(\ADCD)$ and let $\sh{S}:= \sh{T}(\ADCD)$. Since $\ADCD$ is transverse, by Proposition~\ref{prop-ringsnoeth} $S$ and $\sh{S}$ are noetherian.

Let $\sh{K}$ be the left ideal 
\[ \bigoplus_{n\geq 0}  (\sh{A}' \sh{D}^{\sigma^s}\cdots \sh{D}^{\sigma^{n-1}} \sh{C}^{\sigma^n} \cap \sh{I}_p^{\sigma^n} )\Lsh_n \]
of $\sh{S}$, and let
$K:= H^0(X, \sh{K})$.  Let  $I:=T_{\geq 1}$.  Then $I$ is a right ideal of $S$, and it follows from Corollary~\ref{cor-easy3} that $T$ and $\I_S^{\ell}(I)$  are equal in large degree.  Our choice of $p$ forces
$I \cap K / IK$ to be infinite-dimensional; by the right-handed version of \cite[Lemma~2.2]{R-idealizer}, $T$ is not  noetherian.  Likewise, if we let $\sh{I}:= \sh{T}_{\geq 1}$ and compare $\sh{I} \cap \sh{K}$ with $\sh{IK}$, then we may apply the right-handed version of Proposition~\ref{prop-ibma-lnoeth} to conclude that $\sh{T}$ is not noetherian. Thus no such $p$ can exist.
\end{proof}

 Recall that a $\kk$-algebra $T$ is {\em strongly right (left) noetherian} if, for any noetherian commutative $\kk$-algebra $C$, the algebra $T \otimes_{\kk} C$ is right (left) noetherian.  
To end the section, we consider when one of the algebras $T= T(\surfD)$ is strongly right or left noetherian. 

\begin{proposition}\label{prop-SRN}
Let $\surfD = \surfdata$ be transverse ADC data, where $\Lsh$ is $\sigma$-ample.  Let $T:= T(\surfD)$.  Then $T$ is strongly left noetherian if and only if $\sh{D} = \struct_X$ and $\sA$ is invertible, and $T$ is strongly right noetherian if and only if $\sh{D} = \sh{C} = \struct_X$.
\end{proposition}
\begin{proof}
The statements are symmetric, so it suffices to prove the first one.  
Suppose that $\sD= \sO_X$ and $\sA$ is invertible.
Let $\Omega$ be the locally principal Weil divisor defined by $\sA$, and let $\Lsh':= \Lsh(-\Omega+ \sigma^{-1}(\Omega))$.   Then $T$ is 
a subalgebra of 
 the twisted homogeneous coordinate ring $B':= B(X, \Lsh', \sigma)$ as follows.  Let $\sh{J}: = \sh{C}  \sh{I}_{\Omega}$ and let  $J := \bigoplus_{n \geq 1} H^0(X, \Lsh'_n \sJ^{\sigma^n})$.  Then $J$ is a  left ideal of $B'$, and 
 $T$ and $\I^{r}_{B'}(J)$ are equal in large degree.  
 By the left-handed version of \cite[Proposition~7.2]{S-idealizer}, $T$ is strongly left noetherian.  

Let $\dot{\surfD}= \geomdata$ be the geometric data associated to $\surfD$.  Suppose now that $\sh{D} \neq \struct_X$ or that $\sA$ is not invertible; that is, that $\Lambda \cup Z \neq \emptyset$.    Let $U \subset X$ be an open affine subset of $ X$ and let $C := \struct_X(U)$.  
Let 
\[M:= \bigoplus_{n \geq 0} \sh{T}_n(\sigma^{-n}(U)).\]
Note that $U$ contains infinitely many points in $\bigcup_{n\geq 0} \{\sigma^n Z \cup \sigma^n{\Lambda}\}$, and that these points are dense in $U$.
The proof of \cite[Lemma~7.14]{S-idealizer} works in our setting to show that $M$ is a finitely generated left $T \otimes_{\kk} C$-module that is not generically flat over $C$.  By \cite[Theorem~0.1]{ASZ1999}, therefore, $T$ is not strongly left noetherian.
\end{proof}

\section{The $\chi$ conditions}\label{CHI}

We now begin to give homological properties of the rings $T(\surfD)$.  In this section, we focus on the Artin-Zhang $\chi$ conditions.  We first recall the relevant definitions.

\begin{defn}\label{def-chi}
Let $R$ be a finitely generated, connected $\NN$-graded $\kk$-algebra, and let  $j \in \NN$.  We say that {\em $R$ satisfies right $\chi_j$} if, for all $i \leq j$ and for all finitely generated graded right $R$-modules $M$, we have 
\[ \dim_{\kk} \underline{\Ext}^i_{\rgr R}(\kk, M) < \infty.\]
We say that {\em $R$ satisfies right $\chi$} if $R$ satisfies right $\chi_j$ for all $j \in \NN$.  
We similarly define {\em left $\chi_j$} and {\em left $\chi$}; we say $R$ {\em satisfies $\chi$} if it satisfies left and right $\chi$.
\end{defn}

The  condition $\chi_1$ is the most important of the $\chi$ conditions:  if a graded ring $T$ satisfies right $\chi_1$, then it may be reconstructed from the category $\rqgr T$ \cite{AZ1994}.  The higher $\chi$ conditions are more mysterious.  However, if a ring satisfies left and right $\chi$, then it is well-behaved in some significant ways; for example, a ring satisfying $\chi$ has a {\em balanced dualizing complex}, by \cite[Theorem~6.3]{VdB1996}, \cite[Theorem~4.2]{YZ1997}, and thus has a noncommutative version of Serre duality.

We will see that left or right maximality of the ADC data $\surfD$ determine the behaviour in particular of $\chi_1$.  To analyze this, we will need to standardize the ADC data defining our rings slightly more.  The issue is that it is possible to have ADC data $\surfD \neq \surfD'$, with $T(\surfD)$ and $T(\surfD')$  equal in large degree, so that $\surfD$ is right maximal but $\surfD'$ is not.   For example, let $p$ be a point of $X$ with a critically dense $\sigma$-orbit.  Let
\begin{align*}
 \sA & :=  \sI_{\sigma^2(p)} \sI_p,  & \sD& := \sI_p, & \sC &:= \sO_X
\end{align*}
and
\begin{align*}
 \sA' := \sD' & :=  \sI_{\sigma^2(p)} ,   & \sC' &:= \sI_{\sigma^2(p) } \sI_{\sigma(p)}.
\end{align*}
Note that $\sC' \subsetneqq \bigr((\sD'(\sD')^{\sigma}:\sA' \bigl) = \sI_{\sigma(p)}$.  
Let $\surfD = (X, \sL, \sigma, \sA, \sD, \sC, 1)$ and let $\surfD':= (X, \sL, \sigma, \sA', \sD', \sC', 2)$.  Then $\surfD$ is right maximal but $\surfD'$ is not, and $T(\surfD)_{\geq 2} = T(\surfD')_{\geq 2}$.

To correct this, we give the following definition.

\begin{defn}\label{def-standard}
 Let $\surfD = \surfdata$ be transverse ADC data such that the cosupport of $\sA$ is 0-dimensional; let  $\dot{\surfD} = (X, \sigma, \Lambda, Z, \Gamma, \emptyset)$ be the associated geometric data.  If $\Gamma \cap \{\sigma^n(Z)\}_{n<0} = \emptyset$, we say that $\surfD$ is {\em right standard}.  If $\Lambda \cap \{\sigma^n(Z)\}_{n>0}= \emptyset$, we say that $\surfD$ is {\em left standard}.
\end{defn}

In the example above, $\surfD$ is not left standard but $\surfD'$ is left standard.

We leave to the reader the proof of the following elementary lemma:
\begin{lemma}\label{lem-standard}
 Let $\surfD=\surfdata$ be transverse ADC data so that the cosupport of $\sA$ is 0-dimensional.  Then there are left standard ADC data $\surfD'$ and right standard ADC data $\surfD''$ so that $\sT(\surfD)$, $\sT(\surfD')$, and $\sT(\surfD'')$ are equal in large degree. \qed
\end{lemma}

Standardizing ADC data is important because of the following lemma.

\begin{lemma}\label{lem-strongmax}
 Let $\surfD = \surfdata$ be transverse ADC data so that the cosupport of $\sA$ is 0-dimensional, and let $\dot{\surfD} = (X, \sigma, \Lambda, Z, \Gamma, \emptyset)$ be the associated geometric data.  Let $\sT := \sT(\surfD)$.

$(1)$.  If $\surfD$ is left standard and right maximal, then for any $m \gg 0$ we have
\[\shHom_X(\sT_n, \sT_{n+m}^{\sigma^{-m}}) = \sT_m^{\sigma^{-m}}\]
for $n > s$.

$(2)$.  If $\surfD$ is right standard and left maximal, then for any $m \gg 0$ we have
\[ \shHom_X(\sT_n, \sT_{n+m}) = \sT_m^{\sigma^n}\]
for $n > s$.
\end{lemma}
\begin{proof}
 By symmetry, it suffices to prove $(1)$.  Let $m > s$ be such that $\sigma^m(\Lambda) \cap \Lambda = \emptyset$, and let $n > s$.  Then
\begin{multline}\label{BOOM}
\shHom_X(\sT_n, \sT_{n+m}^{\sigma^{-m}}) = \\
\shHom_X(\sA\sD^{\sigma^s}\cdots \sD^{\sigma^{n-1}}\sC^{\sigma^n},
\sA^{\sigma^{-m}} \sD^{\sigma^{-m+s}}\cdots\sD^{\sigma^{n-1}} \sC^{\sigma^n})\sL_m^{\sigma^{-m}}= \\
\shHom_X(\sA, \sA^{\sigma^{-m}} \sD^{\sigma^{-m+s}} \cdots \sD^{\sigma^{s-1}})\sL_m^{\sigma^{-m}}.
\end{multline}
By assumption on $m$ and by left standardness,
\[\Lambda \cap \bigl( \sigma^{-m}(\Lambda) \cup \sigma^{m-s}(Z) \cup \cdots \cup \sigma(Z) \bigr) = \emptyset.\]
So \eqref{BOOM} is equal to 
\[
 \sA^{\sigma^{-m}} \sD^{\sigma^{-m+s}}\cdots \sD^{\sigma^{-1}} (\sD \cdots \sD^{\sigma^{s-1}}: \sA) = \sT_m^{\sigma^{-m}}
\]
by right maximality of $\surfD$.
\end{proof}

We give two more preliminary lemmas.

\begin{lemma}\label{lem-hom}
Let $R$ be a commutative regular local ring of dimension 2 with residue field $\kk$, and let $I \subset J$ be cofinite-dimensional ideals of $R$.  

$(1)$  The natural map 
\[ \Hom_R(J, R) \to \Hom_R(I, R)\]
is an isomorphism.  In particular, 
\[ \Hom_R(J, R) \cong \Hom_R(I, R) \cong \Hom_R(R, R) = R.\]

$(2)$ If $\phi,\psi \in \Hom_R(J, R)$, then $\phi = \psi $ if and only if $\phi|_I = \psi|_I$.
\end{lemma}
\begin{proof}
Since  $R$ is regular of dimension $> 1$, $\Ext^1_R(\kk, R) = 0$.  $(1)$ follows directly, and $(2)$ follows from 
$(1)$.
\end{proof}

Recall that following \cite{Gabriel1962}, if $S$ is a graded ring and $M \in \rGr S$, we  also consider  $M$ to be an element of  $\rQgr S$; we make similar conventions for bimodule algebras.   

\begin{lemma}\label{lem-biff}
Let $\ADCD = \ADCdata$  be transverse ADC data, where $\sh{D} \neq \struct_X$.  Let $\sh{S} := \sh{T}(\ADCD)$ and let $ \sh{B} := \sh{B}(X, \Lsh, \sigma)$.  Then $\Hom_{\rQgr \sh{S}}(\sh{S}, \sh{B}/\sh{S})$ is infinite-dimensional.
\end{lemma}
\begin{proof}
Without loss of generality we may suppose that $\Lsh = \struct_X$.  
Let $t \geq s$ be such that  for $n \geq t$ the ideal sheaves $\sh{A}$ and $\sh{D}^{\sigma^n}$ are comaximal, and $\sh{C}^{\sigma^n}$ and $\sh{D}$ are comaximal.   
Let $n \geq 2t$ 
and let $\sh{E}_n:= \sh{D}^{\sigma^{t}}\cdots \sh{D}^{\sigma^{n-t}}$.  Note that $\sh{E}_n \supseteq \sh{S}_{n+m}$ for any $m \geq 0$.  Since $\struct_X/\sh{E}_n$ and $\struct_X/\sh{S}_{n+m}$ are 0-dimensional and agree at all points in the cosupport of $\sh{E}_n$, the natural map
\[ \struct_X/\sh{S}_{n+m} \to \struct_X/\sh{E}_n\]
splits.  Working pointwise, we obtain maps
\[  \alpha_m: \shHom_X(\sh{S}_n,\struct_X/\sh{E}_n) \to \shHom_X(\sh{S}_{n+m}, \struct_X/\sh{S}_{n+m}).\]

The multiplication maps 
\[ (\struct_X/\sh{S}_n )\otimes \sh{S}_m^{\sigma^n} \to \struct_X/\sh{S}_{n+m}\]
and
\[ \sh{S}_n \otimes \sh{S}_m^{\sigma^n} \to \sh{S}_{n+m}\]
induce maps
\[ \gamma_m:  \shHom_X(\sh{S}_n, \struct_X/\sh{S}_n) 
 \to \shHom_X(\sh{S}_n \otimes \sh{S}_m^{\sigma^n}, \struct_X/\sh{S}_{n+m})\]
and
\[ \delta_m:  \shHom_X(\sh{S}_{n+m}, \struct_X/\sh{S}_{n+m}) \to \shHom_X(\sh{S}_n \otimes \sh{S}_m^{\sigma^n}, \struct_X/\sh{S}_{n+m}).\]
We obtain a diagram
\[ \xymatrix{
\shHom_X(\sh{S}_n, \struct_X/\sh{E}_n ) \ar[d]_{\alpha_m} \ar[r]^{\alpha_0} 
& \shHom_X(\sh{S}_n, \struct_X/\sh{S}_n) \ar[d]^{\gamma_m} \\
\shHom_X(\sh{S}_{n+m}, \struct_X/\sh{S}_{n+m}) \ar[r]_{\delta_m} 
& \shHom_X(\sh{S}_n \otimes \sh{S}_m^{\sigma^n}, \struct_X/\sh{S}_{n+m}) }\]
which is easily seen to be commutative.  Thus taking global sections and summing over $m$,  we obtain, for any $n' \geq n$, an injection
\[ \alpha:   \Hom_X (\sh{S}_n, \struct_X/\sh{E}_n) \to \Hom_{\rGr \sh{S}}(\sh{S}_{\geq n'}, \sh{B}/\sh{S}).\]
Taking the direct limit, we obtain a map
\beq\label{FOO}
 \Hom_X(\sh{S}_n, \struct_X/\sh{E}_n)\to \Hom_{\rQgr \sh{S}}( \sh{S},  \sh{B}/\sh{S}).
\eeq
If $f \in \Hom_X(\sS_n, \sO_X/\sE_n)$ is nonzero, then working pointwise we see that $\alpha_m(f) \in \Hom_X(\sS_{n+m}, \sO_X/\sS_{n+m})$ is also nonzero for any $m$.  Thus the induced element of $\Hom_{\rQgr \sS}( \sS,  \sB/\sS)$ is nonzero; that is, \eqref{FOO} is injective.   
  Thus
\[ \dim_{\kk}\Hom_{\rQgr \sh{S}}(\sh{S}, \sh{B}/\sh{S}) \geq \dim_{\kk} \Hom_X(\sh{S}_n, \struct_X/\sh{E}_n) \geq n-2t.\]
\end{proof}

The next two results describe when various $\chi$ conditions hold for the algebras $T(\surfD)$.

\begin{theorem}\label{thm-chi1}
Let $\surfD = \surfdata$ be  transverse ADC data, where $\Lsh$ is $\sigma$-ample, and let $T := T(\surfD)$.   Let $\dot{\surfD} = \geomdata$ be the associated geometric data.

$(1)$  If $T$ satisfies left $\chi_1$, then  $\surfD$ is left maximal. If $T$ satisfies right $\chi_1$, then  $\surfD$ is right maximal and the cosupport of $\sA$ is 0-dimensional.

$(2)$  If $\surfD$ is maximal, then $T$ satisfies left and right $\chi_1$.

$(3)$  Suppose that  $\surfD$ is left maximal and right standard.    If 
 $\{ \sigma^n \Gamma \}_{n \geq 0}$ is critically transverse, then $T$ satisfies left $\chi_1$. 

$(4)$  Suppose that  $\surfD$ is right maximal and left standard, and that the cosupport of $\sA$ is 0-dimensional.    If 
 $\{ \sigma^n \Lambda \}_{n \leq 0}$ is  critically transverse, then $T$ satisfies right $\chi_1$.
\end{theorem}

\begin{proof}
$(1)$  By symmetry it suffices to prove the first statement.  Suppose that $\surfD$ is not left maximal.  Let $\sh{A}':= (\sh{D}\cdots \sD^{\sigma^{s-1}}: \sh{C})$, let $\mathbb{E} := (X, \Lsh, \sigma, \sh{A}', \sh{D}, \sh{C},s)$ and let $S := T(\mathbb{E})$.  Then $I:= T_{\geq 1}$ is a non-irrelevant right ideal of $S$.  Since $T$ and $\I^{\ell}_S(I)$ 
are equal in large degree and $I \cdot (S/T) = 0$, the factor $S/T$ is an infinite-dimensional torsion left $T$-module.  Therefore,
\[ S/T \hookrightarrow \underline{\Ext}^1_{T \lgr}(\kk, T),\]
and  $T$ fails left $\chi_1$.  

If $\surfD$ is maximal, then it is automatically left and right standard.  Thus $(2)$  follows from $(3)$ and $(4)$.  As $(3)$ and $(4)$ are symmetric, it suffices to prove $(4)$.  

Arguing as in the proof of \cite[Theorem~7.1]{KRS}, to show that $T$ satisfies right $\chi_1$ it suffices to prove that for any coherent right $\sh{T}$-module $\sh{N}$,  the natural map
\beq \label{rainier}
 H^0(X, \sh{N}) \to \bigoplus_m \Hom_{\rqgr \sh{T}}( \sh{T}, \sh{N}[m])
\eeq
has a right bounded cokernel.  Further, since $\sh{N}$ has an ascending chain of submodules whose factors are either Goldie torsion or free, we may assume that either $\sh{N}$ is Goldie torsion or that $\sh{N} = \sh{T}$.  

First suppose that $\sh{N} \neq 0$ is Goldie torsion.  Clearly each $\sh{N}_n$ is a torsion sheaf; since $\sh{N}$ is coherent, generated in degree $\leq n_1$ for some $n_1$, there is a proper subscheme $Y$ of $X$ so that $\Supp \sh{N}_n \subseteq Y$ for all $ n \in \NN$.
  By right maximality, $Z \cup \sigma^{-1}(Z) \cup \cdots \cup \sigma^{-(s-1)}(Z) \supseteq \Gamma$.  Since $\{\sigma^n \Lambda\}_{n \in \ZZ}$ and $\{\sigma^n Z\}_{n \in \ZZ}$ are critically transverse by transversality of $\surfD$ and by assumption on $\Lambda$, there  is $\ell \geq n_1$ so that $Y$ does not meet $\sigma^n Z$ or $\sigma^n \Lambda$ for $n \leq -\ell$.  Thus
$\sh{T}_m^{\sigma^n}|_Y \cong \Lsh_m^{\sigma^n}|_Y$ for all $n \geq \ell$ and $m \geq 0$, and
\[ \sh{T}_n^{\sigma^i}|_Y \cong (\sh{T}^{\sigma^i}_{\ell -i} \sh{T}_{n+i-\ell}^{\sigma^{\ell}})|_Y \cong ( \sh{T}_{\ell-i}^{\sigma^i} \Lsh_{n+i-\ell}^{\sigma^{\ell}})|_Y\]
for all $n \geq \ell$ and $0 \leq i \leq \ell$.  If $i \geq \ell$, then $\sh{T}_n^{\sigma^i}|_Y \cong \Lsh_n^{\sigma^i}|_Y$.    

Let $n \geq \ell$.  Then 
\[
\sh{N}_n = \sum_{i=0}^{n_1} \sh{N}_i \sh{T}_{n-i}^{\sigma^i} 
= \sum_{i=0}^{n_1} \sh{N}_i \sh{T}_{\ell-i}^{\sigma^i} \Lsh_{n-\ell}^{\sigma^{\ell}} 
= \sh{N}_{\ell} \Lsh_{\ell}^{-1} \Lsh_n.\]
 In particular, if $m \geq \ell$ then 
\[ \sh{N}[m]_n = (\sh{N}_{\ell} \Lsh_{\ell}^{-1})^{\sigma^{-m}} \otimes \Lsh_m^{\sigma^{-m}} \otimes \Lsh_n
 \cong \sh{N}_{\ell}^{\sigma^{-m}} \otimes \Lsh_{n+m-\ell}^{\sigma^{\ell-m}}
\]
for all $n \geq 0$, by Lemma~\ref{lem-shift}.

By choice of $\ell$, for any  $m \geq \ell$ we have
\[ \sigma^m(Y) \cap \Bigl(\{ \sigma^n Z\}_{n \leq 0} \cup \{ \sigma^n \Lambda \}_{n \leq 0}\Bigr) = \emptyset.\]
 Then for any $y \in \sigma^m(Y)$ and $n  \geq 0$, we have
\[  (\sh{T}_{n} )_y = (\Lsh_{n})_y.\]
Thus there are maps (in fact isomorphisms)  
\[ \Hom_X(\sh{T}_n, \sh{N}[m]_{n}) \to \Hom_X(\sh{T}_{n+k}, \sh{N}[m]_{n+k})\]
that  induce a map
\beq\label{tiger}
 \Hom_X (\sh{T}_n, \sh{N}[m]_n) \to \Hom_{\rgr \sh{T}}(\sh{T}_{\geq n}, \sh{N}[m])
\eeq
for any $n \geq 0$.  
This is  the inverse of the natural map  
$\Hom_{\rgr \sh{T}}(\sh{T}_{\geq n}, \sh{N}[m]) \to \Hom_X(\sh{T}_n, \sh{N}[m]_n)$, and so \eqref{tiger} is an isomorphism for any $n \geq 0$.  This isomorphism is clearly compatible with the maps in the direct system
\[ \lim_{n \to \infty} \Hom_{\rgr \sh{T}} (\sh{T}_{\geq n}, \sh{N}[m])\]
and so 
\begin{multline*}
 \Hom_{\rqgr \sh{T}}( \sh{T},  \sh{N}[m]) = \lim_{n \to \infty} \Hom_{\rgr \sh{T}} (\sh{T}_{\geq n}, \sh{N}[m]) \\
\cong \Hom_X(\sh{T}_0, \sh{N}[m]_0) \cong H^0(X, \sh{N}_m) = H^0(X, \sT_m^{\sigma^{-m}}).
\end{multline*}
Taking direct limits, we see that \eqref{rainier} is an isomorphism in  degree $\geq \ell$.

Now suppose that $\sh{N} = \sh{T}$.
Using Lemma~\ref{lem-strongmax}(1), choose $m_0 \geq 0$ so that if $m \geq m_0$ and $n \geq s$, then 
$\shHom_X(\sT_n, \sT_{n+m}^{\sigma^{-m}}) = \sT_m^{\sigma^{-m}}$.  If $m \geq m_0$, then $\sh{T}[m]_n = \sh{T}_{m+n}^{\sigma^{-m}}$ for any $n \geq 0$ by Lemma~\ref{lem-shift}.

Fix $n \geq s$ and consider the natural maps
\begin{multline}\label{sthelens} 
\xymatrix {H^0(X, \sh{T}_m) \ar[r]^{\sigma^{-m}}
& H^0(X, \sh{T}[m]_0) \ar[r] &} \\
\xymatrix{ 
& \Hom_{\rgr \sh{T}}(\sh{T}, \sh{T}[m]) \ar[r]
&  \Hom_X(\sh{T}_n, \sh{T}[m]_n) = \Hom_X(\sT_n, \sT_{n+m}^{\sigma^{-m}}). }
\end{multline}
For  $ m \gg 0$, we have by Lemma~\ref{lem-strongmax} that
\[ \Hom_X(\sT_n,\sT_{n+m}^{\sigma^{-m}}) = H^0(X, \sT_m^{\sigma^{-m}}) = H^0(X, \sT[m]_0),\]
and \eqref{sthelens} is an isomorphism.  Thus \eqref{rainier} is an isomorphism in large degree. 
\end{proof}

\begin{theorem}\label{thm-chi2}
Let $\surfD=\surfdata$ be transverse ADC data, where $\Lsh$ is $\sigma$-ample, and let $T:=T(\surfD)$.  The following are equivalent
\begin{enumerate}
\item[$(a)$]
$T$ satisfies right $\chi_2$;
\item[$(b)$]
 $T$ satisfies left $\chi_2$;
\item[$(c)$]
$T$ satisfies $\chi$;
\item[$(d)$]
  $T$ is a twisted homogeneous coordinate ring; that is, $\sh{A} = \sh{D}=\sh{C} = \struct_X$.
\end{enumerate}
\end{theorem}
\begin{proof}
  We show $(a) \iff (c) \iff (d)$; the other implications follow by symmetry.  It is trivial that $(c) \implies (a)$, and $(d) \implies (c)$ follows by \cite[Theorem~6.3]{VdB1997} (or alternately, \cite[Theorem 4.2]{YZ1997}) from the fact, proved in \cite[Theorem~7.3]{Y1992}, that twisted homogeneous coordinate rings have balanced dualizing complexes.  Thus it suffices to prove that $(a) \implies (d)$.

Suppose then that $(d)$ fails, and $T$ is not a twisted homogeneous coordinate ring.   Note that if  $\surfD$ is not right maximal, then $T$ fails right $\chi_2$ by Theorem~\ref{thm-chi1}(1).  
If  $\surfD$ is right maximal and $\sh{D} = \struct_X$, then as $T$ is not a twisted homogeneous coordinate ring we have $\sh{A} \neq \struct_X$.  In this case
 $T$  fails right $\chi_2$ by  \cite[Proposition~8.4(2)]{S-idealizer}.  Thus  it suffices to suppose that $\sh{D}$ is  nontrivial and $\surfD$ is right maximal, and show that $T$ fails right $\chi_2$.  

Let $B: =B(X, \Lsh, \sigma) $ and let $\sh{B}:=\sh{B}(X, \Lsh, \sigma)$.
We first claim that 
\beq \label{cc1}
\Hom_{\rQgr T}(T, B) \cong 
\Hom_{\rQgr \sh{T}}(\sh{T}, \sh{B}) \cong \kk.
\eeq

The first isomorphism is a consequence of the equivalence of categories in Theorem~\ref{thm-VdBSerre}.  Thus it suffices to prove the second.  We may without loss of generality suppose that $\Lsh = \struct_X$.  Note that  $\Hom_X(\sh{T}_n, \struct_X) \cong \Hom_X(\struct_X, \struct_X) = \kk$ for all $n$.  

Fix $n \geq 0$, and let $\phi: \sh{T}_{\geq n} \to \sh{B}$ be a right $\sh{T}$-module homomorphism.  We claim that $\phi$ is determined by $\phi|_{\sh{T}_n}$.  So suppose  that $\psi: \sh{T}_{\geq n} \to \sh{B}$ is another right $\sh{T}$-module map.  For all $i \geq n$, let $\phi_i = \phi|_{\sh{T}_i}$, and similarly for $\psi$.  Suppose that $\psi_n = \phi_n$, and let $i \geq 1$.  Consider the maps
\[ \xymatrix{
\sh{T}_n \otimes \sh{T}_i^{\sigma^n} \ar[r]^(.6){\phi_n \otimes 1}& \sh{T}_{i}^{\sigma^n} \ar[r]^{\alpha} & \struct_X, 
}\]
where $\alpha$ is the canonical inclusion.  This factors as
\[ \xymatrix{
\sh{T}_n \otimes \sh{T}_i^{\sigma^n} \ar[r]^(.6){\phi_n \otimes 1} \ar[rd]_{\beta} 
& \sh{T}_i^{\sigma^n} \ar[r]^{\alpha}
&  \struct_X \\ 
& \sh{T}_n \cdot \sh{T}_i^{\sigma^n} \ar[ru]_{\gamma} \ar[u]
}\]
where $\beta$ is the canonical map of $\struct_X$-modules.  Note that $\gamma$ is simply 
$\phi_n|_{\sh{T}_n \sh{T}_i^{\sigma^n}}$.  Furthermore, as $\phi$ is a right $\sh{T}$-module map, we have that
\[ \gamma = \phi_{n+i}|_{\sh{T}_n \cdot \sh{T}_i^{\sigma^n} \subset \sh{T}_{n+i}}.\]
Repeating this analysis for $\psi$, we see that 
\[ \psi_{n+i}|_{\sh{T}_n\sh{T}_i^{\sigma^n}} = \phi_{n+i}|_{\sh{T}_n\sh{T}_i^{\sigma^n}}.\]

Let $\dot{\surfD} = \geomdata$ be the geometric data associated to $\surfD$.  
Transversality of $\surfD$ implies that all points in $\Lambda \cup Z \cup \Gamma$ have dense orbits; in particular, they are contained in the smooth locus of $X$.  
By Lemma~\ref{lem-hom}, therefore, $\psi_{n+i} = \phi_{n+i}$.

The canonical map 
\[ \Hom_{\rGr \sh{T}}(\sh{T}_{\geq n}, \sh{B}) \to \Hom_X(\sh{T}_n, \struct_X)\]
is therefore injective.  Since $\Hom_{\rGr \sh{T}}(\sh{T}_{\geq n}, \sh{B})\neq 0$, we have 
\beq\label{bar}
\Hom_{\rGr \sh{T}}(\sh{T}_{\geq n}, \sh{B}) \cong \kk
\eeq
for any $n \geq 0$.  A similar argument shows that the diagram
\beq\label{foo}
 \xymatrix{
\Hom_{\rGr \sh{T}}(\sh{T}_{\geq n}, \sh{B})  \ar[rr] \ar[d]_{\cong} &&
\Hom_{\rGr \sh{T}}(\sh{T}_{\geq n+1}, \sh{B}) \ar[d]^{\cong} \\
\Hom_X(\sh{T}_n, \struct_X)  && \Hom_X(\sh{T}_{n+1}, \struct_X) \\
&\Hom_X(\struct_X, \struct_X)\ar[ru]^{\cong} \ar[lu]_{\cong} &
 }\eeq
commutes.  In particular, the top row of \eqref{foo} is an isomorphism.  

Now, 
\[ \Hom_{\rQgr \sh{T}}(\sh{T}, \sh{B}) \cong \lim_{n \to \infty} \Hom_{\rGr \sh{T}}(\sh{T}_{\geq n}, \sh{B}).\]
The maps in the direct system are precisely those in the top row of \eqref{foo}, and so they are all isomorphisms.  By \eqref{bar}, we have that $\Hom_{\rQgr \sh{T}}(\sh{T}, \sh{B}) \cong \kk$.

By Lemma~\ref{lem-biff}, $\Hom_{\rQgr \sh{T}}(\sh{T}, \sh{B}/\sh{T})$ is infinite-dimensional. 
 From the long exact $\Hom$ sequence 
\[ \Hom_{\rQgr \sh{T}}( \sh{T}, \sh{B}) \to \Hom_{\rQgr \sh{T}}(\sh{T}, \sh{B}/\sh{T}) \to \Ext^1_{\rQgr \sh{T}}(\sh{T}, \sh{T})\]
we deduce  that
\[\Ext^1_{\rQgr \sh{T}}(\sh{T}, \sh{T}) \cong \Ext^1_{\rQgr T}(T, T)\]
is infinite-dimensional.  By \cite[$(\dag)$, p. 274]{AZ1994}, $T$ fails right $\chi_2$.
\end{proof}

We believe that  maximal ADC rings are the proper generalizations of na\"ive blowups at a point, even though they may not be generated in degree 1.   The poor behavior of non-maximal ADC rings is evidence for this opinion.  
In  \cite[Example~5.1]{RS-0} Rogalski and Stafford construct a  na\"ive blowup algebra  that satisfies $\chi_1$ on the right but not on the left.  In our terms, this example may also be constructed as follows.  Let $X = \PP^2$, and let $\sigma \in \Aut(X)$ be such that the point $p$ has a critically dense orbit.  Let $(x,y)$ be local coordinates at $p$.  Let $\Lsh$ be any ample (and therefore $\sigma$-ample) invertible sheaf on $X$.  We define three ideal sheaves cosupported at $p$.  Let $\sh{A}$ be defined by 
\[ \sh{A}_p := (x,y) = \sI_p.\]
Let $\sh{D} := \sh{A}^3$, so $\sh{D}_p = (x,y)^3$.  Let
\[ \sh{C}_p := (x^2, y^2).\]
Then the transverse ADC data $\ADCD = (X, \sL, \sigma, \sA, \sD, \sC, 1)$  is right standard and left but not right maximal, as $\sC \subsetneqq (\sD:\sA) = \sA^2$.  By Theorem~\ref{thm-chi1}, the ring $S:= S(\ADCD)$ satisfies left $\chi_1$ and fails right $\chi_1$.  Note that $\sh{S}:=\sh{T}(\ADCD)$ satisfies
\[ \sh{S}_n = \sh{S}_1 \sh{S}_1^{\sigma} \cdots \sh{S}_1^{\sigma^{n-1}}\]
and that for sufficiently ample $\Lsh$ the ring $T$ is generated in degree 1.  From the perspective of the current paper, the surprisingly pathological properties of some  na\"ive blowups noted in \cite{RS-0} thus come from the non-maximality of the associated ADC data.  Note that, by \cite[Theorem~1.1]{KRS} (or by the results of this section), a na\"ive blowup at a point on a critically dense orbit always satisfies left and right $\chi_1$ and fails left and right $\chi_2$.

\section{Noncommutative projective schemes}\label{PROJ}
Let $\surfD =\surfdata$ be transverse ADC data, and let $T:= T(\surfD)$.  In this section, we determine the homological properties of the category $\rQgr T$ (and, by symmetry, $T \lQgr$).  As pointed out by Artin and Zhang \cite{AZ1994}, this category, or, more properly, the pair $(\rqgr T,  T)$ is the correct noncommutative analogue of Proj of a finitely generated commutative graded ring.  We are particularly interested in studying what we informally call the {\em cohomological dimension} of the category; that is, the cohomological dimension of the ``global sections'' functor 
\[ \Hom_{\rQgr T}( T, \blank).\]
We show that this dimension is finite for the rings $T(\surfD)$.

We begin the section by showing that the category $\rqgr T$ depends only on $X, \sigma$, and $\sh{D}$.  In particular, there is a na\"ive blowup  (or twisted homogeneous coordinate ring) $S$, at a scheme cosupported on points with infinite distinct orbits, with $\rqgr T \simeq \rqgr S$.

We recall a result of Rogalski on idealizers in graded algebras.

\begin{proposition}\label{prop-R-catequiv}
{\em \cite[Lemma~3.2]{R-idealizer}}
Let $U$ be a noetherian connected  $\NN$-graded $\kk$-algebra, let $H$ be a graded left ideal of $U$ so that $\dim_{\kk}  (U/H) = \infty$, and let $V := \I^{r}_U(H)$.  
Assume in addition that ${}_V U$ is finitely generated and that $\dim_{\kk}(V/H) < \infty$.  
Then the functor
\[ {}_U(H \otimes_V \blank): V \lqgr \to U \lqgr\]
is an equivalence of categories with quasi-inverse
\[ {\rm res}: {}_U M \mapsto {}_V M.\]
Further, the functor
\[ (\blank \otimes_U H)_V:  \rqgr U \to \rqgr V\]
is an equivalence of categories, with quasi-inverse 
\[ (\blank \otimes_V U): \rqgr V \to \rqgr U.\]
\qed
\end{proposition}

Let $\surfD$ be transverse ADC data.  Using Proposition~\ref{prop-R-catequiv}, we construct a na\"ive blowup $S$ so that $\rqgr S \simeq \rqgr T$.

\begin{theorem}\label{thm-catequiv}
Let $\surfD = \surfdata$ be  transverse ADC data, where $\Lsh$ is $\sigma$-ample.  Let $\sh{T} := \sh{T}(\surfD)$ and let $T := T(\surfD)$.  Let $Z$ be the subscheme defined by $\sh{D}$ and let $\sh{S}$ be the na\"ive blowup $\sh{S} := \sh{S}(X, \Lsh, \sigma, Z)$.  Let $S := S(X, \Lsh, \sigma, Z)$.  Then the categories $\rqgr \sh{T}$, $\rqgr T$, $\rqgr S$, and $\rqgr \sh{S}$ are equivalent.  Likewise, the categories $\sh{T} \lqgr$, $T \lqgr$, $S \lqgr$, and $\sh{S} \lqgr$ are equivalent.  
\end{theorem}

Note that $T$ may not be an idealizer in $S$.   However, we may still obtain this result from repeated applications of Proposition~\ref{prop-R-catequiv}.  
\begin{proof}
By Theorem~\ref{thm-VdBSerre} and by symmetry, it suffices to prove that $\rqgr T \simeq \rqgr S$.

We first note that $T$ is an idealizer in a maximal ADC  ring.  
Let $\sh{C}' := (\sh{D} : \sh{A})$ and let
$\sh{A}':= (\sh{D}\cdots\sD^{\sigma^{s-1}} : \sh{C}')$.  
Let $\mathbb{E} := (X, \Lsh, \sigma, \sh{A}, \sh{D}, \sh{C}',s)$ and let
$\mathbb{F}:=(X, \Lsh, \sigma, \sh{A}', \sh{D}, \sh{C}',s)$.  
Note that  $\mb{F}$ is maximal.  
Let $R := T(\mathbb{E})$ and let $U := T(\mathbb{F})$.  Let $I := T_{\geq 1}$ and let $J := R_{\geq 1}$.  Then
$I$ is a left ideal of $R$ and $J$ is a  right ideal of $U$.

By an easy generalization of Corollary~\ref{cor-easy3}, $T$ and $\I^{r}_R(I)$ 
are equal in large degree, and 
$R$ and $\I^{\ell}_U(J)$ 
are equal in large degree.  Applying Proposition~\ref{prop-R-catequiv} twice, 
we obtain equivalences of categories 
\[ \rqgr T \simeq \rqgr U\]
and
\[ T \lqgr \simeq U \lqgr.\]

Now, $\sh{A}' \supseteq \sh{D}\cdots \sD^{\sigma^{s-1}}$.  Consider the transverse ADC data
\[ \mb{G} := (X, \Lsh, \sigma, \sD, \sD, \sC',1).\] 
Let $V:= T(\mb{G})$.  
Let $K := V_{\geq 1}$.  Then $K$ is a right ideal of $U$ and a left ideal of $S = S(X, \Lsh, \sigma, Z)$.  Further, 
\[ V= \I^{\ell}_U(K) = \I^{r}_S(K),\] 
as a consequence of our assumptions on $\sh{A}, \sh{D}, \sh{C}$.  Therefore, applying Proposition~\ref{prop-R-catequiv} again, we see that
\[ \rqgr U \simeq \rqgr S\]
and
\[ U \lqgr \simeq S \lqgr.\]
\end{proof}

We note that the equivalences from Theorem~\ref{thm-catequiv} do not take the distinguished object $ T \in \rqgr T$ to 
$ S\in \rqgr S$.  Unpacking the functors from Proposition~\ref{prop-R-catequiv}, the equivalence
$\rqgr T \to \rqgr S$ is given by:
\[ \xymatrix@R=0pt{
\rqgr T  \ar[r] & \rqgr R \ar[r]&  \rqgr U \ar[r] &  \rqgr V \ar[r]& \rqgr S \\
 M_T \ar@{|->}[r] &  M\otimes_T R_R \ar@{|->}[r] & M \otimes_T J_U \ar@{|->}[r]  & M \otimes_T J_V \ar@{|->}[r] & M \otimes_T J  \otimes_V S_S, } \]
and the equivalence $T \lqgr \to S \lqgr$ is given by
\[ \xymatrix@R=0pt@C=8pt{
T \lqgr  \ar[r] & R \lqgr \ar[r] & U \lqgr \ar[r] & V \lqgr \ar[r] & S \lqgr \\
{}_TN \ar@{|->}[r] & {}_RI \otimes_T N \ar@{|->}[r] & {}_UU \otimes_R I \otimes_T N \ar@{|->}[r] & {}_VK \otimes_R I \otimes_T N \ar@{|->}[r] & {}_SK \otimes_V K \otimes_R I \otimes_T N. }\]

\begin{corollary}\label{cor-cat-invariants}
Let $\surfD =\surfdata$ be transverse ADC data, where $\Lsh$ is $\sigma$-ample.  Let $T := T(\surfD)$.  The categories $\rqgr T$ and $T \lqgr$ depend only on $X, \sigma$, and $\sh{D}$.
\end{corollary}
\begin{proof}
Let $Z$ be the subscheme defined by $\sh{D}$.  
Let $\sh{S} := \sh{S}(X, \Lsh, \sigma, Z)$ be the na\"ive blowup bimodule algebra at $Z$.  Since $\rqgr T \simeq \rqgr \sh{S}$, the category $\rqgr T$ depends only on $X, \Lsh, \sigma$, and $\sh{D}$ (or $Z$).  By Lemma~\ref{lem-catequiv}, however, $\rgr \sh{S}$ does not depend on $\Lsh$; thus neither does $\rqgr \sh{S}$.
\end{proof}

This proves part (2) of Theorem~\ref{ithm-qgr}.

One consequence of Theorem~\ref{thm-catequiv} is that if $R$ is an ADC ring, then $\rqgr R$ is equivalent to $\rqgr S$ for some na\"ive blowup $S$.  Since it is shown in  \cite[Theorem~6.7]{KRS} that the category of Goldie torsion modules over $S$ is equivalent (in Proj) to the category of torsion quasi-coherent sheaves on $X$, it follows that the Goldie torsion subcategory of $\rqgr R$ will be equivalent to the category of torsion sheaves on $X$.  We record this directly as:

\begin{proposition}\label{prop-torsion}
Let $\surfD = \ADCdata$ be  transverse maximal ADC data.  Let $\sh{S} := \sh{T}(\surfD)$.  Then the functor
\[ F: \sh{N} \mapsto \bigoplus_{n \geq 0} \sh{N} \otimes \Lsh_n\]
from $\struct_X \lMod \to \rQgr \sh{S}$
restricts to  an equivalence between the 
full subcategory $GT(\struct_X \lmod)$ of coherent  torsion sheaves on $X$ and the full subcategory
$GT(\rqgr \sh{S})$ of objects in $\rqgr \sh{S}$ that are images of Goldie torsion right $\sh{S}$-modules.
\end{proposition}
\begin{proof}
We essentially follow the proof of \cite[Theorem~6.7]{KRS}, even though $\sh{S}$ is not generated in degree 1.  Without loss of generality, we may assume that $\Lsh = \struct_X$.  Let $\sh{F} = \bigoplus_{n \geq 0} (\sh{F}_n)_{\sigma^n}$ be a nonzero coherent graded Goldie torsion $\sh{S}$-module; we may assume that $\sF$ is torsion-free as an $\sS$-module.   Then $\sh{F}$ is generated in degree $\leq n_1$ for some $n_1$.  Since the $\sh{F}_n$ are clearly  torsion sheaves on $X$,  there is some proper subscheme $Y$ of $X$ so that $\Supp \sh{F}_n \subseteq Y$ for all $n \in \NN$.    By transversality, let $n_2 \geq n_1$ be such that 
\[ \sh{S}_m^{\sigma^n}|_Y \cong \struct_Y\]
for all $n \geq n_2$ and all $m \geq 0$.  This means that
\[
 \sh{F}_n \otimes \sh{S}^{\sigma^n}_{m-n}  \cong \sh{F}_n \otimes \struct_X \]
for all $m \geq n \geq n_2$.  Therefore, we have
\[ \sh{F}_j \otimes \sh{S}_{n_2-j}^{\sigma^j} \otimes \sh{S}_{n-n_2}^{\sigma^{n_2}} \cong
\sh{F}_j \otimes \sh{S}_{n_2-j}^{\sigma^j} \cdot \sh{S}_{n-n_2}^{\sigma^{n_2}} \cong
\sh{F}_j \otimes \sh{S}_{n_2-j}^{\sigma^j}\cdot \struct_X\]
for all $n \geq n_2 \geq j$.  This implies that for $n \geq n_2$ we have $\sh{F}_n = \sh{F}_{n_2}$.  We may apply the proof of \cite[Theorem~6.7]{KRS} to our situation.  Just as in that proof, it follows that $F$ takes coherent Goldie torsion to coherent objects, is surjective on Goldie torsion objects, and is full and faithful on morphisms between Goldie torsion objects.   Thus $F$ restricts to an equivalence, as claimed.
\end{proof}

We now begin to investigate the cohomological and homological dimensions of the category $\rqgr T$.  
\begin{proposition}\label{prop-homdim}
Let $\surfD = \surfdata$ be transverse ADC data, where $\Lsh$ is $\sigma$-ample and $X$ is smooth.  Then both $\rQgr T$ and $T \lQgr$ have finite homological dimension:  that is, there is some $i$ so that for $j >i$, we have
\[ \Ext^j_{\rQgr T}(\blank, \blank) = 0,\]
and similarly for $T \lQgr$.
\end{proposition}
\begin{proof}
It suffices to prove the statement for $\rqgr T$.  If $T$ is a na\"ive blowup, this is \cite[Theorem~6.8]{RS-0}.  For general $T$, the result follows from Theorem~\ref{thm-catequiv}.
\end{proof}

In contrast to homological dimension, the cohomological dimension of the functor $\Hom_{\rQgr T}(T, \blank)$ depends on the distinguished object $ T$, and thus  is not necessarily preserved under the category equivalences from Theorem~\ref{thm-catequiv}.  Stafford and Van den Bergh asked \cite[page 194]{SV} whether any noetherian ring must have finite cohomological dimension.  In \cite[Example~9.7]{S-idealizer}, we gave an example of a right, but not left, noetherian ring $R$ so that the right cohomological dimension of $R$ is infinite.  The ring $R$ is a geometric idealizer  defined by non-transverse  data on a singular surface.  

We conclude this section by showing that, in contrast, for  algebras of transverse ADC data, left and right cohomological dimension are always finite.  Before giving this result, we prove a vanishing lemma for a certain class of $\Ext$ groups.
  
\begin{lemma}\label{lem-torsion}
Let $\mb{E} = \ADCdata$ be  transverse  maximal ADC data, where $\Lsh$ is $\sigma$-ample.  Let $Y \subset X$ be a proper subscheme such that $Y$ is locally principal at every singular point of $X$.    

 Let 
$U := T(\mb{E})$.  Let $\sh{K}$ be a quasicoherent torsion sheaf on $X$, and let 
\[ K := \bigoplus_{n \geq 0} H^0(X, \sh{K} \otimes \Lsh_n).\]
Let 
\[ J' := \bigoplus_{n \geq 0} H^0(X, \sh{I}_Y \Lsh_n)\]
and let \[ J := J'\cap U.\]
Then 
$ \Ext^i_{\rQgr U}(J, K) = 0 $
for all $i \geq 5$.
\end{lemma}
\begin{proof}
Let 
\[ \sh{U} := \sh{T}(\mb{E}).\]
Let 
\[ \sh{J} := \bigoplus_{n \geq 0} (\sh{I}_Y \Lsh_n \cap \sh{U}_n)\]
and let $\sh{B} := \sh{B}(X, \Lsh, \sigma)$.  
By  Theorem~\ref{thm-VdBSerre}, it suffices to prove that
\[ \Ext^i_{\rQgr \sh{U}}(\sh{J}, \sh{K} \otimes \sh{B}) = 0\]
for $i \geq 5$.  To show this, without loss of generality we may assume that $\Lsh = \struct_X$.

We first suppose that $\sh{K}$ is in fact supported on $X^{\sing}$.     Let $\sh{H}$ be the reflexive hull $\sh{I}_Y^{**}$.  
Our assumption on $Y$ implies that $\sh{H}$ is invertible.
We have $\sJ_n \subseteq \sI_Y \subseteq \sH$, and an induced exact sequence
\beq\label{seq-seq} 
 \shHom_X(\sH/\sJ_n, \sK) \to \shHom_X(\sH, \sK) \to \shHom_X(\sJ_n, \sK) \to \shExt_X^1(\sH/\sJ_n, \sK).
\eeq
The cosupport of any $\sh{U}_n$ is disjoint from $X^{\sing}$.  Thus $\sH/\sJ_n$ is supported on a finite set disjoint from $X^{\sing}$.  The first and last terms of \eqref{seq-seq} are therefore 0, and we have
\[ \Hom_X(\sJ_n, \sK)  \cong \Hom_X(\sH, \sK)\]
for any $n$.  
   There is clearly an induced isomorphism 
\[ \xymatrix{
\Hom_X(\sh{H}, \sh{K}) \ar[r]^(.4){\cong} & \Hom_{\rGr \sh{U}}(\sh{J}_{\geq n}, \sh{K} \otimes \sh{B}).}\]
The inverse of this map is the canonical restriction 
\[ \Hom_{\rGr \sh{U}}(\sh{J}_{\geq n}, \sh{K} \otimes \sh{B}) \to \Hom_X(\sh{J}_n, \sh{K}) \cong \Hom_X(\sh{H}, \sh{K}).\]
Since this isomorphism exists for any $n \geq 0$, we see that
\[ \Hom_{\rQgr \sh{U}}(\sh{J}, \sh{K} \otimes \sh{B}) \cong
 \lim_{n \to \infty}\Hom_{\rGr \sh{U}}(\sh{J}_{\geq n}, \sh{K} \otimes \sh{B}) \cong 
\Hom_X(\sh{H}, \sh{K}).\]
If $\sh{K} \to \sh{K}'$ is a morphism of sheaves, where
$\sh{K}'$ is another quasicoherent sheaf supported on 
$X^{\sing}$, then the reader may check that the diagram
\[ \xymatrix{
\Hom_{\rQgr \sh{U}}(\sh{J}, \sh{K} \otimes \sh{B}) \ar[r] \ar[d]
& \Hom_{\rQgr \sh{U}}(\sh{J}, \sh{K}' \otimes \sh{B}) \ar[d] \\
\Hom_X(\sh{H}, \sh{K}) \ar[r] & \Hom_X(\sh{H}, \sh{K}') }\]
commutes.   Therefore, the two functors
\[ \Hom_{\rQgr \sh{U}}(\sh{J}, \blank \otimes \sh{B})\]
and
\[ \Hom_X(\sh{H}, \blank)\]
are isomorphic as functors from the category of sheaves supported on $X^{\sing}$ to $\rm{Ab}$.

Let 
\[ \sh{K} \to \sh{I}_{\bullet}\]
be a minimal injective resolution of $\sh{K}$; note that each term of $\sh{I}_{\bullet}$ is also supported on $X^{\sing}$ and in particular is torsion.  Then the cohomology of 
\beq \label{swine}
 \Hom_{\rQgr \sh{U}}(\sh{J}, \sh{I}_{\bullet} \otimes \sh{B}) \cong \Hom_X(\sh{H}, \sh{I}_{\bullet})
\eeq
computes  the groups
\[ \Ext^i_X(\sh{H}, \sh{K}).\]
On the other hand, $\sh{B}$ is a flat $\struct_X$-module.  Thus $\sh{I}_{\bullet} \otimes \sh{B}$ is a resolution of $\sh{K} \otimes \sh{B}$ as a $\sh{B}$-module and therefore as a $\sh{U}$-module.  Proposition~\ref{prop-torsion} implies that it is an injective resolution in $\GT(\rQgr \sh{U})$.  Since an object in $\GT(\rQgr \sh{U})$ is injective if and only if it is injective as an object of $\rQgr \sh{U}$, the cohomology of \eqref{swine}
also computes the groups
\[ \Ext^i_{\rQgr \sh{U}}(\sh{J}, \sh{K} \otimes \sh{B}).\]
Now, if $i \geq 2$ then we have
\[ \Ext^i_X(\sh{H}, \sh{K}) \cong H^i(X, \sh{K} \otimes \sh{H}^{-1}) = 0\]
as the support of $\sh{K}$ has dimension at most 1.  Thus
\[ \Ext^i_{\rQgr \sh{U}}(\sh{J}, \sh{K} \otimes \sh{B}) = 0\]
if $i \geq 2$, and certainly the lemma holds for $\sh{K}$.  

Now let $\sh{K}$ be any quasicoherent torsion sheaf, and let 
\[ 0 \to \sh{K} \to \sh{I}_0 \to \sh{I}_1 \to \sh{I}_2 \to \cdots \]
be a minimal injective resolution of $\sh{K}$.  Since minimal injective resolutions commute with localization, the sheaves $\sh{I}_n$ for $n\geq 3$ are supported on $X^{\sing}$.  Let $\sh{K}'$ be the cokernel of the map $\sh{I}_1 \to \sh{I}_2$.  Then $\sh{K}' \hookrightarrow \sh{I}_3$ is supported on $X^{\sing}$. 
By Proposition~\ref{prop-torsion},
\[ 0 \to \sh{K}\otimes\sh{B} \to \sh{I}_0 \otimes \sh{B} \to \sh{I}_1 \otimes \sh{B} \to \cdots \]
is an injective resolution of $\sh{K}\otimes \sh{B}$ in $\rQgr \sh{U}$.   Thus
\[ \Ext^i_{\rQgr \sh{U}}(\sh{J}, \sh{K} \otimes \sh{B}) \cong \Ext^{i-3}_{\rQgr \sh{U}}(\sh{J}, \sh{K}'\otimes \sh{B})\]
for $i \geq 4$.   We have seen that this vanishes if $i-3 \geq 2$. 
\end{proof}

\begin{proposition}\label{prop-cohdim}
Let $\surfD = \surfdata$ be transverse ADC data, where $\Lsh$ is $\sigma$-ample.  Let $T := T(\mb{D})$.  Then $T$ has finite left and right cohomological dimension.
\end{proposition}
\begin{proof}
By symmetry, it suffices to prove the statement on the right.
We show there is some $i$ so that
\[ 0 = \Ext^i_{\rQgr T}( T,  M) = \lim_{n \to \infty} \Ext^i_{\rGr T}(T_{\geq n}, M)\]
for any $M \in \rgr T$.

As in the proof of Theorem~\ref{thm-catequiv}, let $\sh{C}':= (\sD \sD^{\sigma} \cdots \sD^{\sigma^{s-1}} : \sh{A})$ and let 
$\sh{A}':= (\sD \sD^{\sigma} \cdots \sD^{\sigma^{s-1}}: \sh{C}')$.  Let 
\[ \mb{F}:=(X, \Lsh, \sigma, \sh{A}', \sh{D}, \sh{C}',s),\]
so $\mb{F}$ is maximal ADC data.
Let $\sh{U}$ be the ADC bimodule algebra $\sh{U}:= \sh{T}(\mb{F})$ and let $U := T(\mb{F})$.  Let $\mb{E}:= (X, \sL, \sigma, \sA, \sD, \sC',s)$, and let $\sR:=\sT(\mb{E})$ and $R:= T(\mb{E})$.   

Note that $\sR_{\geq 1}$ is a right ideal of $\sU$ and a finitely generated left $\sT$-module, since ${}_{\sT} \sR$ is finitely generated.
By Corollary~\ref{cor-easy3}, $\sR = \I^{\ell}_{\sU}(\sR_{\geq 1})$ 
in large degree, and $\sT = \I^{r}_{\sR}(\sT_{\geq 1})$ 
in large degree.
Similarly, $R = \I^{\ell}_U(R_{\geq 1})$ 
in large degree, and $T = \I^{r}_R(T_{\geq 1})$ 
in large degree.

Let $J:= R_{\geq 1}$, so $J$ is a right $U$-module and a left $T$-module.  
Recall from the comments after the proof of Theorem~\ref{thm-catequiv} that the functor
\[ M_T \mapsto M \otimes_T J_U\]
induces an equivalence of categories between $\rqgr T$ and $\rqgr U$.  In particular,
\[ \Ext^i_{\rQgr T} (T, M) \cong \Ext^i_{\rQgr U}(J, M \otimes_T J).\]
Thus it suffices to prove that for $i \gg 0$, we have $\Ext^i_{\rQgr U}(J, N) = 0$ for any $N_U$.   Further, since $E(N)/N$ is Goldie torsion, it is enough to prove this for $N$ Goldie torsion.  This follows directly from Proposition~\ref{prop-torsion} and the previous lemma, since the subscheme  defined by $\sA$ is locally principal at every point of $X^{\sing}$.
\end{proof}

We conjecture that if $\surfD$ is transverse, then the correct value for the cohomological dimension of $\rQgr T$ is 2.

\section{Maximal orders}\label{MAXORD}

In this section, we study ADC rings of maximal transverse data.  We show that ADC rings on normal surfaces are maximal orders:  the noncommutative version of an integrally closed ring.  This is a new class of maximal order, not previously observed.

\begin{theorem}\label{thm-maxord}
Let $\mb{E}= (X, \sL, \sigma, \sA, \sD, \sC, 1)$ be  transverse maximal ADC data, where $\Lsh$ is $\sigma$-ample, and 
further suppose that $X$ is normal.  Then $T := T(\mb{E})$ is a maximal order.
\end{theorem}

For example, let $X$ be a normal surface, let $\sigma \in \Aut(X)$, let $\sL$ be a $\sigma$-ample invertible sheaf on $X$, and let $p \in X$ have a critically dense orbit.  Let $\sA = \sC = \sI_p$.  Let $x,y$ be local coordinates at $p$,  and let $\sD$ be the ideal sheaf cosupported at $p$ so that $\sD_p = (x,y^2)$.  Then 
$\mb{E} = (X, \sL, \sigma, \sA, \sD, \sC, 1)$ is maximal ADC data, and so $T(\mb{E})$ is a maximal order.  Since $\sA\sC \subsetneqq \sD$, no Veronese of $T(\mb{E})$ is generated in degree 1.   By \cite[Proposition~3.18]{RS-0}, $T(\mb{E})$ is not a na\"ive blowup algebra.

  We will work inside the graded quotient ring $D:= Q_{gr}(T)$ of $T$; note that $D \cong K[z,z^{-1}; \sigma]$ where $K$ is the function field of $X$.  In this section, we will use $z$ as a dummy variable to indicate degree.  That is, we let $\sh{T}:= \sh{T}(\ADCD)$, and we write:
\[ T = \bigoplus_{n \geq 0} H^0(X, \sh{T}_n) z^n.\]
The advantage of this notation is that we now have a natural inclusion 
\[ T \subset K[z, z^{-1}; \sigma]\]
and so we may write
\[ D = Q_{gr}(T) = K[z, z^{-1};\sigma].\]
Note that this convention requires us to write
\[ B(X, \Lsh, \sigma) = \bigoplus_{n\geq 0} H^0(X, \Lsh_n) z^n.\]

We consider now what happens when $s$ is arbitrary.

\begin{proposition}\label{prop-maxord}
Let $\ADCD = \ADCdata$ be transverse ADC data, where $\Lsh$ is $\sigma$-ample, $X$ is normal, and the cosupport of $\sA$ is 0-dimensional.  
Define 
\[ \sh{I}_n := \sh{A} \sh{D}^{\sigma^s} \cdots \sh{D}^{\sigma^{n-1}} \sh{C}^{\sigma^n}\]
for all $n \geq s$, and let 
\[ S := \kk + \bigoplus_{n \geq s} H^0(X, \sh{I}_n \Lsh_n) z^n.\]
Let $K$ be the function field of $X$ and let $D := K[z, z^{-1}; \sigma]$ be the graded quotient ring of $S$.  
Then the following are equivalent:

$(1)$  $S \ver{s}$ is a maximal order;

$(2)$  $\ADCD$ is maximal ADC data;

$(3)$ $S$ has finite codimension in a graded maximal order.  That is, there is a graded ring $D \supset R \supseteq S$, with $R/S$ finite-dimensional, so that $R$ is a maximal  order.
\end{proposition}

Before giving the proof, we prove a preparatory lemma.  Recall that if $R$ is a noetherian domain with (full) quotient ring $Q$, and $J$ is an ideal of $R$, we define
\[ O_{\ell}(J) :=\{q \in Q \st qJ \subseteq J\} \]
and 
\[ O_r(J) := \{q \in Q \st Jq \subseteq J\}.\]
By \cite[Lemma~3.1.12]{MR}, $O_{\ell}(J)$ and $O_r(J)$ are equivalent orders to $R$.   Further, by \cite[Proposition~5.1.4]{MR}, $R$ is a maximal order if and only if  $O_{\ell}(J) = O_r(J) =  R$ for all ideals $J$ of $R$.  

\begin{lemma}\label{lem-overring}
Let $\ADCD = \ADCdata$ be transverse  ADC data, where $X$ is a normal surface, $\Lsh$ is $\sigma$-ample, and the cosupport of $\sA$ is 0-dimensional.  Let  $B: =B(X, \Lsh, \sigma)$.  
Let $T:=T(\ADCD)$ and let $S \subseteq B$ be a graded subring with $S_{\geq t} = T_{\geq t}$ for some $t \geq 0$.  
Let $D:= K[z, z^{-1}; \sigma]$ be the graded quotient ring of $B$. Let $R$ be a graded overring of $S$ so that
\[ B_{< t} + S \supseteq R \supseteq S.\]
Let $J$ be a graded ideal of $R$.   Then $O_{\ell}(J)$ and $O_r(J)$ are graded subrings of $B$. 

 If $\ADCD$ is maximal, then  
\[ B_{<s}+T \supseteq O_{\ell}(J) + O_r(J).\]
\end{lemma}
\begin{proof}
It suffices to prove the result for $O_{\ell}(J)$.  Let $Q$ be the full quotient ring of $R$;  we may naturally embed $D$ in $Q$.  The proof of \cite[Lemma~9.1]{R-generic} shows that $O_{\ell}(J) \subseteq D$.  It is obviously graded.

Let $\sh{T}:=\sh{T}(\ADCD)$.
Let $Z$ be the cosupport of $\sh{D}$.  By Proposition~\ref{prop-ideals},  there are a $\sigma$-invariant ideal sheaf $\sh{J}$ on $X$, cosupported away from orbits of points in $Z$, and an integer $k$ so that
\[ J_n = H^0(X, \sh{JT}_n) z^n\]
for any $n \geq k$.

Let $m \geq k, t, s$ be such that $\sh{JT}_n$ is globally generated for $n \geq m$.  
Now, 
\[
 (O_{\ell}(J)_n  z^{-n}) \cdot (J_m z^{-m})^{\sigma^n} \subseteq J_{n+m}z^{-n-m} \subseteq K.
\]
Multiplying by $\sO_X$ and using the fact that $\sJ\sT_{n}$ and $\sJ\sT_{n+m}$ are globally generated, we obtain
\[ O_{\ell}(J)_n z^{-n} \subseteq \Hom_X(\sh{JT}_{m}^{\sigma^n}, \sh{JT}_{n+m}) \]
for any $n \in \NN$.  Since $X$ is normal,
\[ \shHom_X(\sh{JT}_{m}^{\sigma^n}, \sh{JT}_{n+m}) = \shHom_X(\sh{T}_{m}^{\sigma^n}, \sh{T}_{n+m}) \subseteq \Lsh_n\]
for any $n$, so $O_{\ell}(J)\subseteq B$.
Further,
\beq\label{sss} 
\shHom_X(\sh{T}_{m}^{\sigma^n}, \sh{T}_{n+m}) = (\sh{AD}^{\sigma^s} \cdots \sh{D}^{\sigma^{n+m-1}} \sh{C}^{\sigma^{n+m}}: \sh{A}^{\sigma^n} \sh{D}^{\sigma^{n+s}} \cdots \sh{D}^{\sigma^{n+m-1}} \sh{C}^{\sigma^{n+m}}) \Lsh_n.
\eeq

Assume now that $\ADCD$ is maximal, and let $n \geq s$.  Then  \eqref{sss} is equal to 
\[  \sh{AD}^{\sigma^s}\cdots\sh{D}^{\sigma^{n-1}}(\sh{D}^{\sigma^n}\sh{D}^{\sigma^{n+1}} \cdots \sh{D}^{\sigma^{n+s-1}}:\sh{A}^{\sigma^n}) \Lsh_n= \sh{AD}^{\sigma^{s}}\cdots \sh{D}^{\sigma^{n-1}} \sh{C}^{\sigma^n}\Lsh_n.\]
  That is,
\eqref{sss} is equal to $ \sh{T}_n$
and $O_{\ell}(J)_{\geq s} = T_{\geq s}$ as claimed.
\end{proof}

\begin{proof}[Proof of Proposition~\ref{prop-maxord} and Theorem~\ref{thm-maxord}]
Note that Theorem~\ref{thm-maxord} is the $s=1$ case of Proposition~\ref{prop-maxord}.

$(3)$ $\Rightarrow$ $(2)$.
We prove the contrapositive.  
Let $D \supset R \supseteq S$ be a graded overring so that $R/S$ is finite-dimensional.  Suppose that $\ADCD$ is non-maximal ADC data; we will show that there is an equivalent order $T \supsetneqq R$.  In fact, we will see that $T/S$ is infinite-dimensional.  

Since $R/S$ is finite-dimensional, we have $R\cdot S_{\geq t} \subseteq S_{\geq t}$ for some $t$.  
By the previous lemma, $R \subseteq O_{\ell}(S_{\geq t}) \subseteq B(X, \Lsh, \sigma)$. 
Let $a \geq 1$ be minimal so that $R_{a'} \neq 0$ for any $a' \geq a$.  Then for any $n \geq a$, define the ideal sheaf $\sh{J}_n$ 
to be the image of the natural map $(R_n z^{-n}) \otimes \Lsh_n^{-1} \to \struct_X$; that is, $\sh{J}_n$ is the base ideal of the sections  in $R_n z^{-n} \subseteq H^0(X, \Lsh_n)$.  

Since $(2)$ fails, there is some point $p$  so that   either
\beq \label{foo1}
 \sh{A}_{p} \subsetneq (\sh{D} \cdots \sh{D}^{s-1}: \sh{C})_{p}
\eeq 
or
\[ 
\sh{C}_{p} \subsetneq (\sh{D} \cdots \sh{D}^{s-1}: \sh{A})_{p}.
\]
The two cases are symmetric; we will suppose that \eqref{foo1} holds and show that $R$ is not a maximal order. 

Let $Z$ be the cosupport of $\sh{D}$. We first suppose that the orbit of $p$ meets $Z$.  
Let $\struct := \struct_{X,p}$.  For any $i \in \ZZ$, define $p_i := \sigma^{-i}(p)$.  For any $j$, we may identify the stalk $\struct_{X, p_j}$ with $\struct$.  
Using this identification, define for any $i \in \NN$ and $j \in \ZZ$ an ideal 
\[ \mf{r}^i_j := (\sh{J}_i)_{p_j} \subseteq \struct.\]
The multiplication law on $R$ translates to the equation
\beq\label{eq-fund}
 \mf{r}^i_j \mf{r}^k_{j-i} \subseteq \mf{r}^{i+k}_j
\eeq
for any $j$ and $i,k \geq 1$.  Let $m \geq s$ be such that for $n \geq m$, the sheaf $\sh{R}_n = \sh{I}_n \Lsh_n$ is globally generated.  
 Then
\[ \mf{r}^n_j = (\sh{I}_n)_{p_j}\]
for any $n \geq m$ and any $j$.  

It follows from \eqref{eq-fund} that there are integers $b \leq 0$ and $c \geq s$ so that $\mf{r}^n_j = \struct$ for any $j \not\in [b, n+c)$.  By reindexing the orbit of $p$, and possibly changing $\sh{A}$ and $\sh{C}$, we may assume that $b = 0$ and $c=s$.  (We leave the tedious but routine verification to the reader.)  Thus $ (\sh{AC})_{p_j} = \struct $
for any $j \not\in [0, s)$.

 For $0 \leq k \leq s-1$, define the following ideals of $\struct$:
\[ \mf{a}_k := \sh{A}_{p_k}\]
and
\[ \mf{c}_k := \sh{C}_{p_k}.\]
Let 
\[ \mf{d} := \sh{D}_p.\]
For $n \geq m$ we have
\[ \mf{r}^n_k = \mf{a}_k \mbox{ if  $n \geq m$ and $ 0 \leq k < s$},\]
\[ \mf{r}^n_k = \mf{d} \mbox{ if $n \geq m$ and $  s \leq k < n$},\]
and
\[ \mf{r}^n_k = \mf{c}_{k-n} \mbox{ if $n \geq m$ and $ n \leq k < n+s$}.\]
If $k \not \in [0, n+s)$ and $n \geq a$, then $\mf{r}^n_k = \struct$.  
Let
\[ \sh{A}' :=  (\sh{D} \cdots \sh{D}^{s-1}: \sh{C}).\]
By assumption, $\sh{A}' \supsetneqq \sh{A}$.  
Define
\[ \mf{a}'_k := (\sh{A}')_{p_k} = (\mf{d}:\mf{c}_k)\]
for $0 \leq k < s$.  

Let 
\[ T:=   R + \bigoplus_{n \geq m} H^0(X, \sh{A}' \sh{D}^{\sigma^s} \cdots \sh{D}^{\sigma^{n-1}} \sh{C}^{\sigma^n}\Lsh_n).\]
Then $(T/R)_n \neq 0$ for all $n \gg 0$, so  certainly $T /S$ is infinite-dimensional.  Since $R_s T \subset R$,  if $T$ is a ring it is  an equivalent order to $R$ (and $S$).

We thus must show that  $T$ is multiplicatively closed.   It suffices to show that
\[ R_i T_j + T_j R_i \subseteq T_{i+j} \]
for $1 \leq i < m$ and $j \geq m$.
For $n \geq m$, let
\[ \mf{s}^n_j := (\sh{A}' \sh{D}^{\sigma^s} \cdots \sh{D}^{\sigma^{n-1}} \sh{C}^{\sigma^n})_{p_j}.\]
That is, 
\begin{align}  \label{biffboff}\begin{split}
\mf{s}^n_j = \mf{a}'_j & \mbox{ for $0 \leq j < s$,} \\
\mf{s}^n_j = \mf{d} & \mbox{ for $s \leq j < n$,  and} \\ 
\mf{s}^n_j = \mf{c}_{j-n} & \mbox{ for $n \leq j < n+s$.}
\end{split} \end{align}
Note also that 
\beq\label{wow}
\mf{s}^n_j = \mf{r}^n_j \mbox{ if $n \geq m$ and $j \not\in [0, s)$}.
\eeq

We must check that
\beq\label{wanted}
 \mf{r}^i_k \mf{s}^j_{k-i} + \mf{s}^j_k \mf{r}^i_{j-k} \subseteq \mf{s}^{i+j}_k
\eeq
for any $1 \leq i < m$, $j \geq m$, and $k \in \ZZ$.

We first show that 
\beq \label{eq1}
\mf{s}^j_k \mf{r}^i_{j-k} \subseteq \mf{s}^{i+j}_k.
\eeq
If $k <j$ then \eqref{biffboff} gives that $\mf{s}^j_k = \mf{s}^{i+j}_k$, so \eqref{eq1} is automatic.  And if $k \geq j$,  from \eqref{wow} we have
$\mf{s}^j_k = \mf{r}^j_k$ and $\mf{s}^{i+j}_k = \mf{r}^{i+j}_k$, and  \eqref{eq1} follows from the fact that $R$ is multiplicatively closed.

We now show that
\beq\label{eq2}
\mf{r}^i_k \mf{s}^j_{k-i} \subseteq \mf{s}^{i+j}_k
\eeq
for any $ 1 \leq i < m$, $j \geq m$, and $k \in \ZZ$.  This is an argument by cases.  If $k < 0$ then $\mf{s}^{i+j}_k = \struct$, so \eqref{eq2} is automatic.  If $0 \leq k \leq s-1$, then  \eqref{wow} gives that
$\mf{s}^m_{m+k} = \mf{r}^m_{m+k}$ and $\mf{r}^{m+i}_{m+k} = \mf{s}^{m+i}_{m+k}$.  Since
 $R$ and $T_{\geq m}$ are multiplicatively closed, we have
\[ \mf{s}^m_{m+k} \mf{r}^i_k \mf{s}^j_{k-i} = \mf{r}^m_{m+k} \mf{r}^i_k \mf{s}^j_{k-i} 
\subseteq \mf{r}^{m+i}_{m+k} \mf{s}^j_{k-i} = \mf{s}^{m+i}_{m+k} \mf{s}^j_{k-i}
\subseteq \mf{s}^{m+i+j}_{m+k}.\]
Now, $k < s < i+j$, so  $\mf{s}^{m+i+j}_{m+k} = \mf{d}$ and $\mf{s}^m_{m+k} = \mf{c}_k$.  Thus we have
\[ \mf{r}^i_k \mf{s}^j_{k-i} \subseteq (\mf{s}^{m+i+j}_{m+k}: \mf{s}^m_{m+k}) = 
(\mf{d}: \mf{c}_k) = \mf{a}'_k .\]
But this is equal to $\mf{s}^{i+j}_k$, so \eqref{eq2} holds.

If $s \leq k < s+i$ then $i+j>k \geq s$, so $\mf{s}^{i+j}_k = \mf{d}$.  We have
\[ \mf{r}^i_k = \struct \cdot \mf{r}^i_k = \mf{r}^m_{m+k} \mf{r}^i_k \subseteq \mf{r}^{m+i}_{m+k} = \mf{s}^{m+i}_{m+k}.\]
Thus
\[ \mf{r}^i_k \mf{s}^j_{k-i} \subseteq \mf{s}^{m+i}_{m+k} \mf{s}^j_{k-i} 
\subseteq \mf{s}^{m+i+j}_{m+k} = \mf{d} = \mf{s}^{i+j}_k,\]
and  \eqref{eq2} holds.  

Finally, if $s+i \leq k$, then note that
\[ \mf{s}^j_{k-i} = \mf{s}^{i+j}_k,\]
so \eqref{eq2} is automatic.  Thus \eqref{eq2} holds in all cases, and $T$ is multiplicatively closed.

This proof also shows that $T\ver{s}/S\ver{s}$ is nonzero and $S\ver{s}$ is an equivalent order to $T\ver{s}$, so  $(1) \Rightarrow (2)$.

$(2) \implies (3)$.  Consider the set of all graded subrings $R$ of $B$ so that  
\[ R_{\geq s} = S_{\geq s}.\]
Since $(B_{< s}+ S)/S$ is finite-dimensional, this set has a maximal element, say $R'$.  Let $J$ be a graded ideal of $R'$.  By Lemma~\ref{lem-overring}, $R' \subseteq O_{\ell}(J) \subseteq B_{<s}+T$.  Maximality of $R'$ therefore implies that $O_{\ell}(J) = R'$, and by symmetry $O_r(J) = R'$.  By \cite[Lemma~9.1]{R-generic}, $R'$ is a maximal order.  

Note that  if $(2)$ holds,  the previous paragraph shows  that $S\ver{s}$ is a maximal order, and $(2) \Rightarrow (1)$.   If $s=1$ and $\ADCD$ is maximal, then $S$ itself is a maximal order.  
\end{proof}

We note that if $\surfD = \surfdata$ where the cosupport of $\sA$ is 1-dimensional, then $T(\surfD)$ is easily seen to be neither a left or right maximal order.

In Section~\ref{CHI}, we commented that maximal ADC algebras are probably the best generalization of na\"ive blowups at a single point, since they satisfy left and right $\chi_1$ automatically.  Theorem~\ref{thm-maxord} gives further proof of this; note that in \cite[Theorem~9.5]{R-generic}, it is shown that a na\"ive blowup of a point in $\PP^2$ is a maximal order.  

There are technical issues that may make it more difficult to work with an ADC algebra  than with a na\"ive blowup at a point, but most of these are relatively easily overcome, as we have seen.  Notably, if $\ADCD$ is maximal ADC data with $\sh{AC} \subsetneqq \sh{D}\cdots \sD^{\sigma^{s-1}}$, then no Veronese of $S(\ADCD)$ is generated in degree 1.  Previously, idealizers were the only observed class of geometric algebras with this property.  The class of ADC algebras thus delimits poor homological properties, such as the failure of $\chi_1$, from failure to be generated in degree 1.

These observations suggest that developing techniques in graded ring theory that do not require the algebras under study to be generated in degree 1 may be important for future research.  This was part of the motivation for the companion paper \cite{S-surfclass}.

\bibliographystyle{amsalpha}
\providecommand{\bysame}{\leavevmode\hbox to3em{\hrulefill}\thinspace}
\providecommand{\MR}{\relax\ifhmode\unskip\space\fi MR }
\providecommand{\MRhref}[2]{%
  \href{http://www.ams.org/mathscinet-getitem?mr=#1}{#2}
}
\providecommand{\href}[2]{#2}

\end{document}